\documentclass[a4paper,11pt]{article}
\usepackage[T1]{fontenc}
\usepackage[utf8]{inputenc}
\usepackage{amsmath,amssymb,amsthm}
\usepackage{booktabs,array,graphicx}
\usepackage[margin=25mm]{geometry}
\usepackage[hidelinks]{hyperref}
\usepackage{float}
\usepackage{url}
\theoremstyle{plain}
\newtheorem{proposition}{Proposition}
\theoremstyle{definition}
\newtheorem{finding}{Finding}
\newtheorem{remark}{Remark}
\newcommand{\abs}[1]{\left|#1\right|}

\newcommand{\OzI}{Ozaki scheme~I}
\newcommand{\Nmod}{N_{\mathrm{mod}}}
\newcommand{\OzII}{Ozaki scheme~II}

\title{Accelerating Multiple-Precision LU Decomposition\\
with Ozaki Scheme II\\
\large --- Multi-component and Arbitrary Precision on CPUs and GPUs ---}
\author{Tomonori Kouya\\
\small Otemon Gakuin University}
\date{\today}

\begin{document}
\maketitle

\begin{abstract}
Solving ill-conditioned linear systems needs LU decomposition in precisions
beyond binary64.  The standard approaches---GMP/MPFR, or multi-component
arithmetic such as double-double---leave every scalar multiply-add inside the
$O(n^3)$ update in multiple precision, and so cannot exploit the low-precision
matrix engines that dominate current hardware.  We build a blocked, partially
pivoted LU decomposition on \OzII{}, which replaces the multiple-precision
GEMM by exact integer modular products followed by an explicit CRT
reconstruction, and implement it for multi-component (DD/TD/QD) and arbitrary
precision, on CPUs and GPUs.  Two contributions make this practical: a direct
conversion between the non-overlapping expansion format and the internal
fixed-point representation, which removes the MPFR round trip and is worth a
factor of $2.4$--$4.7$; and new FP16, FP8 and binary64 GPU back-ends, the FP8
one using a balanced base-$17$ two-digit encoding that keeps the full
$362.8$-bit CRT capacity of INT8.  On an Arm/GB10 and an x86/H100 the fastest
back-end changes with the machine: INT8 wins on GB10, whereas on H100
binary64 is fastest for QD, beating both the native implementation and INT8.
The essential point of \OzII{} is thus not to use a low-precision engine but
to choose the format maximising bits-per-modulus times engine throughput.
For the Lotkin matrix at $p\approx1.2\log_2\mathrm{cond}(A)$ we reach relative
errors of $10^{-645}$ at $n=2048$, up to $2.84\times$ faster than a fully
OpenMP-parallel multiple-precision LU.  We also show by measurement that the
$O(p)$ advantage of scheme II over scheme I applies only to the GEMM term:
modular reduction and CRT grow as $O(p^2)$ and dominate the runtime at the
matrix sizes considered here.  The implementation is released as open source.
\end{abstract}

\tableofcontents

\section{Introduction}
\label{sec:intro}

Solving ill-conditioned linear systems requires an LU decomposition in
multiple precision, beyond the 53 bits of binary64.  The standard way to
realise this is either a variable-precision library such as
GMP/MPFR~\cite{mpfr2007,gmp}, or multi-component (double-double and similar)
arithmetic of the QD~\cite{qd} family, parallelised with OpenMP.  In both
cases, however, every individual scalar multiply-add inside the $O(n^3)$
update remains a multiple-precision operation, so none of the overwhelming
throughput of low-precision matrix engines---AMX-INT8 on CPUs, INT8/FP4 tensor
cores on GPUs---can be exploited.

\subsection{Prior work: Strassen-type algorithms and their limits}
\label{sec:strassen}

The author has previously investigated Strassen's algorithm and the Winograd
variant as a route to faster multiple-precision matrix
multiplication~\cite{kouya2014,kouya2016}.  This classical route, which
reduces the operation count from $O(n^{3})$ to $O(n^{\log_2 7})$, tends to be
relatively favourable in multiple precision because a single multiplication is
expensive, and it did deliver a measurable speedup.  In the arbitrary-precision
regime of hundreds to thousands of bits, however, the gain fell short of
expectations, for two reasons.

\begin{itemize}
\item Strassen's algorithm trades multiplications for additions and
  subtractions.  In multiple precision these are also expensive
  (proportional to the significand length), and the ratio of an addition to a
  multiplication does not shrink as the significand grows.  The benefit of
  reducing the operation count is therefore blunted.
\item Every level of recursive splitting allocates and copies temporary
  matrices, so the volume of memory traffic grows in proportion to the
  significand length.  At arbitrary precision a single element can reach
  hundreds of bytes, and this $O(n^2)$ data movement is no longer negligible.
\end{itemize}

Strassen's algorithm also carries a numerical-stability penalty: in our
measurements the error grows by hundreds to thousands of ulps relative to a
naive implementation.  In short, the Strassen family reduces the
\emph{number} of operations but leaves the \emph{cost of one operation}
untouched, and this is a fundamental limitation in the arbitrary-precision
regime.

By contrast, the family of methods due to Ozaki et
al.~\cite{ozaki2012,ozaki2-2025} decomposes the inputs into low-precision
quantities and replaces the multiple-precision arithmetic itself by several
calls to a highly optimised low-precision GEMM.  This attacks the cost of one
operation directly, and is thus complementary to the Strassen family.  In
\cite{kouya2023} the author reported an optimisation of LU decomposition based
on \OzI{} (splitting the significand into several slices), but \OzI{} needs a
number of products that grows quadratically with the target precision
($O(p^2)$), which is unfavourable at high precision.  \OzII{}, the subject of
this paper~\cite{ozaki2-2025}, uses integer modular decomposition and the
Chinese remainder theorem (CRT) to keep the number of products linear in the
target precision ($O(p)$).  This advantage is however confined to the number
of products, that is, to the term proportional to $n^3$: the pre- and
post-processing for modular reduction and CRT, which is proportional to $n^2$,
grows as $O(p^2)$.  Which of the two dominates the wall-clock time depends on
the combination of $n$ and $p$ (Finding~\ref{find:phase}).

\paragraph{Premise: the stagnation of binary64 performance}
The starting point of this work is a hardware trend of the past decade or so:
binary64 arithmetic performance has been growing slowly relative to that of
low-precision engines.  Kashi et al.~\cite{kashi2024} survey the throughput of
data-centre GPUs across both the NVIDIA and AMD lines and show that, since the
introduction of tensor cores, only INT8/FP16/BF16/TF32 performance has grown
disproportionately; they conclude that ``if the performance growth trends for
low-precision and double-precision arithmetic continue to diverge\dots it may
eventually be possible, or even necessary, to directly approximate
double-precision-accurate basic operations''.  Indeed, the FP64 peak of NVIDIA
data-centre GPUs rose from $19.5/39$~TFLOP/s (vector/tensor) on the A100
(2020) to $33.5/67$~TFLOP/s on the H100 (2022), but then flattened or slightly
declined to $30$~TFLOP/s on the B100 (2024)~\cite{morgan2025}.
Matsuoka~\cite{matsuoka2026} points out that in the Blackwell Ultra generation
native FP64 falls to roughly $1.3$~TFLOP/s while FP8 tensor throughput reaches
the PFLOP/s range, and argues that FP64 should be demoted from a hardware
requirement to an accuracy guarantee obtained by composition---specifically,
by \OzII{}.

Our own measurements agree with this trend.  Swapping only the back-end within
the same \OzII{} code, the effective GEMM throughput on the Arm/GB10 CPU is
$4661$~Gop/s for INT8 against $208$~Gop/s for binary64, a factor of $22.4$
(Table~\ref{tab:phase}).  On the GPU side the effect is even more extreme: on
Arm/GB10 the FP64 units are weak enough that the same LU decomposition runs
$13.0\times$ slower on the GPU than on the machine's own CPU
(Finding~\ref{find:gb10_gpu_slow}).

Work on emulating binary64 with low-precision engines falls into two lines.
One puts \OzI{} on integer matrix-multiply units (INT8): ozIMMU by Ootomo et
al.~\cite{ozimmu2024} and the performance-improved version by Uchino et
al.~\cite{uchino2025}.  The other is the \OzII{} line treated here, in which
Uchino et al.~\cite{ozfp8-2026} achieve binary64-equivalent emulation through
FP8 quantisation.  Note in particular that \cite{ozimmu2024,uchino2025} are
based on \OzI{}, not on \OzII{}.  All of these target binary64-equivalent
accuracy, whereas this paper addresses the high-precision regime from DD/TD/QD
up to several thousand bits.  As shown later (\S\ref{sec:fp816}), this
difference in target precision changes which low-precision engine is optimal.

\subsection{Contributions}

\begin{itemize}
\item Building on a GEMM library based on \OzII{}, we implemented a blocked LU
  decomposition with partial pivoting on both CPU and GPU
  (\S\ref{sec:alg}).
\item We support both multi-component precision (DD/TD/QD) and arbitrary
  precision (MPFR, \texttt{cu\_freal}), and compare CPU and GPU at identical
  precision and identical size (\S\ref{sec:bench}).
\item Using the ill-conditioned Lotkin matrix, we show on both CPU and GPU
  that arbitrary-precision LU decomposition can be solved accurately and
  quickly (\S\ref{sec:lotkin}).
\item In addition to INT8, we implemented new FP16, FP8 and binary64 back-ends
  on the GPU, verified that each yields the same number of correct bits as
  INT8, and compared their performance (\S\ref{sec:fp816},
  \S\ref{sec:fp64}).  FP8 secures the same $362.8$-bit CRT capacity as INT8
  through the balanced base-$17$ two-digit decomposition of
  Eq.~\eqref{eq:fp8digits}.
\item From this comparison we show that the fastest back-end changes with the
  machine.  INT8 is fastest on Arm/GB10, but on x86/H100, whose FP64 units are
  strong, the binary64 back-end is fastest for QD and beats both the existing
  implementation and INT8 (Finding~\ref{find:fp64}).  The essential point of
  \OzII{} is therefore not ``use a low-precision engine'' but ``choose the
  format that maximises the product of bits-per-modulus and engine
  throughput''.  We also give a quantitative reason for not adopting FP32 or
  TF32 (\S\ref{sec:fp32}).
\end{itemize}

\section{Target multiple-precision libraries}
\label{sec:libs}

We describe the four libraries that our implementation handles as input,
output, or internal representation.

\subsection{dtq (CPU, multi-component precision)}
dtq~\cite{dtq} is a derivative of the QD library~\cite{qd} providing
\texttt{dd\_real} (DD, about 106 bits), \texttt{td\_real} (TD, about 159 bits)
and \texttt{qd\_real} (QD, about 212 bits), formed by combining 2, 3 or 4
binary64 numbers.  The internal representation is a non-overlapping expansion:
the leading component is the rounded binary64 value and each subsequent
component is the residual, rounded to binary64 in turn.  Because the
operations are implemented as expanded sequences of binary64 arithmetic with
no memory indirection, they are considerably faster than MPFR at comparable
precision.  A C API operating on raw \texttt{double} arrays is available in
\texttt{c\_dd.h}/\texttt{c\_td.h}/\texttt{c\_qd.h} and can be used directly
from our library, which is written in C.  Since \texttt{dd\_real} and the
others hold their components in a contiguous \texttt{double x[nc]}, such an
array can be viewed directly as an array of \texttt{double[nc]}.  We exploit
this property to implement a path that converts the expansion format directly
into the internal fixed-point representation of \OzII{}
(\S\ref{sec:directconv}, below).

\subsection{MPFR (CPU, arbitrary precision)}
MPFR~\cite{mpfr2007} is a variable-precision floating-point library built on
GMP~\cite{gmp}; the significand length can be set arbitrarily at run time and
correct rounding is guaranteed for every operation.  Our CPU-side \OzII{} GEMM
and LU decomposition take \texttt{mpfr\_t} arrays as input and output.  Being
able to change the significand length at run time makes it suitable for
choosing the precision according to the condition number, as in
\S\ref{sec:lotkin}.

\subsection{gdtq (GPU, multi-component precision)}
gdtq~\cite{gdtq} is the GPU counterpart of dtq, providing \texttt{gdd\_real},
\texttt{gtd\_real} and \texttt{gqd\_real}, which are respectively
\texttt{double2}, \texttt{double3} and \texttt{double4}.  The memory layout of
these types is identical to dtq's \texttt{x[2]}, \texttt{x[3]} and
\texttt{x[4]}, so values can be passed between CPU and GPU without
conversion.  We found on the actual hardware that gdtq's addition and
subtraction are available as \texttt{\_\_host\_\_ \_\_device\_\_}, but the
multiplication kernel (\texttt{two\_prod}) is \texttt{\_\_device\_\_}-only and
cannot be called from the host.  Our \OzII{} path handles the components of
the expansion directly as integers (\S\ref{sec:directconv}, below) and
therefore never calls gdtq arithmetic, so it is unaffected by this
restriction.  The gdtq-native LU decomposition used as a comparison target, on
the other hand, uses \texttt{operator*} and the like inside device kernels,
for which \texttt{\_\_device\_\_}-only is not a problem.

\subsection{MPC\_CUDA (GPU, fixed precision)}
\texttt{cu\_freal<PB>} of MPC\_CUDA~\cite{mpccuda} is a register-resident fixed-precision
floating-point type on the GPU, with $PB$ a compile-time constant ($PB$ a
multiple of 32, $N=\lceil PB/64\rceil$ limbs).  It rounds bit-compatibly with
MPFR's RNDN and is header-only.  Our GPU-side LU decomposition uses it as the
element type; $PB=128/160/224$ correspond to DD/TD/QD precision
respectively.  The conversion to and from MPFR was newly implemented in this
work, and a bit-exact round trip was verified over 2000 trials for
$PB\in\{128,256,512,1024\}$.

\section{\OzII{} and its application to LU decomposition}
\label{sec:alg}

We apply \OzII{} to GEMM and to LU decomposition, taking the four libraries of
the previous section as input and output.  Sections~\ref{sec:gemm}
through~\ref{sec:lu} describe the algorithm independently of the
representation format; \S\ref{sec:impl} then gives the mapping onto each
library and each arithmetic engine.

We fix notation first.  $p$ is the target precision in bits; $t$ (=$Q$) is
the width of the internal fixed-point significand; $w$ is the width of one
chunk of that significand; $S$ is the number of chunks (slices), so $Q=wS$.
$\Nmod(w)$ is the number of moduli required for one pass of width $w$, with
$m_1,\ldots,m_{\Nmod}$ the individual moduli and
$\mathcal P=\prod_\ell m_\ell$ their product.  Matrices are $M\times N\times
K$ for GEMM ($K$ the inner dimension) and $n\times n$ for LU; for square GEMM
we write $M=N=K=n$.  $b$ is the LU panel width.  Note that the triangular
factors $L,U$ of the LU decomposition and the slice count $S$ are unrelated.

\subsection{\OzII{} for GEMM}
\label{sec:gemm}

Given multiple-precision matrices $A\in\mathbb{F}^{m\times k}$ and
$B\in\mathbb{F}^{k\times n}$ at a target precision of $p$ bits, the product
$C=AB$ is obtained with a single rounding to $p$ bits as follows.

\begin{enumerate}
\item[(S1)] Choose shared exponents $\mu_i,\nu_j$ per row and per column.
\item[(S2)] Convert to $Q$-bit signed fixed-point integers as
  $\bar a_{ih}=\mathrm{trunc}(2^{\mu_i}a_{ih})$ and
  $\bar b_{hj}=\mathrm{trunc}(2^{\nu_j}b_{hj})$.
\item[(S3)] Form the centered residues $\bar A^{(\ell)},\bar B^{(\ell)}$ for
  each of the pairwise coprime small moduli $m_1,\ldots,m_{\Nmod}$ ($\le255$).
\item[(S4)] Compute $C^{(\ell)}=\bar A^{(\ell)}\bar B^{(\ell)}$ with an exact
  low-precision GEMM, INT8$\times$INT8$\to$INT32.
\item[(S5)] Reconstruct $\bar c_{ij}$ exactly by CRT.
\item[(S6)] Only now round:
  $c_{ij}=\mathrm{fl}_p(\bar c_{ij}\cdot2^{-(\mu_i+\nu_j)})$.
\end{enumerate}

The integer products in (S4) are exact, and (S6) is the only rounding in the
entire inner product.  This is the source of the accuracy advantage over a
naive multiple-precision GEMM, which rounds at each of the $k$ multiply-adds.

\subsection{Arbitrary precision: splitting into several slices and CRT}
\label{sec:crt}

The uniqueness condition $2\abs{\bar A}\abs{\bar B}<\mathcal P$ for the
modulus product $\mathcal P=\prod_\ell m_\ell$ bounds the fixed-point width
$Q$ that a single CRT pass can reach.  At most 54 pairwise coprime moduli not
exceeding $255$ (the largest power of each prime) are available, giving
$\log_2\mathcal P\approx362.8$ bits in total.

\begin{proposition}[Precision limit of a single pass]
\label{prop:pmax}
With inner dimension $K$ and $g$ guard bits,
\begin{equation}
Q_{\max}=\frac{362.8-\log_2 2K-2}{2},\qquad
p_{\max}\approx Q_{\max}-\log_2 2K-g,
\end{equation}
so that for $K=4096$ and $g=8$ we get only $p_{\max}\approx155$ bits.
\end{proposition}

By Proposition~\ref{prop:pmax}, DD-equivalent precision (106 bits) barely fits
into a single pass, but TD/QD (159/212 bits) and above do not.  \OzII{}
resolves this by splitting the fixed-point value into limbs (several slices).

\paragraph{Slice splitting}
Write the target fixed-point width as $Q=w\cdot S$ and split the $Q$-bit
integers $X,Y$ into $S$ balanced (signed) chunks of $w$ bits:
\begin{equation}
X=\sum_{t=0}^{S-1}v_t\,2^{wt},\qquad \abs{v_t}\le2^{w-1}.
\end{equation}
The product is grouped by digit group $g=t_a+t_b$,
\begin{equation}
XY=\sum_{g=0}^{2S-2}2^{wg}\!\!\sum_{t_a+t_b=g}\!\!v_{t_a}u_{t_b}.
\end{equation}
The lower digit groups ($g<S-1$) are truncated and only the $S$ groups
$g=S-1,\ldots,2S-2$ are retained.  The number of pairs satisfying
$t_a+t_b=g$ is $g+1$ for $g\le S-1$ and $2S-1-g$ for $g\ge S-1$, so the total
number of retained terms is
\begin{equation}
\sum_{g=S-1}^{2S-2}(2S-1-g)=\sum_{i=0}^{S-1}(S-i)=\frac{S(S+1)}{2}.
\end{equation}
Proposition~\ref{prop:trunc} justifies this truncation.

\begin{proposition}[Truncation error of the lower digit groups]
\label{prop:trunc}
For $\abs{v_t},\abs{u_s}\le2^{w-1}$ the truncated part is bounded by
\begin{equation}
\Bigl|\sum_{g=0}^{S-2}2^{wg}S_g\Bigr|
\;\le\;2^{2w-2}\sum_{g=0}^{S-2}(g+1)2^{wg}
\;\le\;(S-1)\,2^{Q-1},
\qquad Q=wL .
\label{eq:truncbound}
\end{equation}
The range of an inner product of inner dimension $K$ is
$\abs{\sum_h X_{ih}Y_{hj}}\le K\,2^{2Q}$, so the relative contribution of the
truncation to that range is at most $(S-1)2^{-Q-1}$.  Consequently
\begin{equation}
Q \;\ge\; p+\log_2\Gamma+\log_2 2K+\log_2 S+g_0 ,
\label{eq:Qcond}
\end{equation}
where $\Gamma$ is the cancellation amplification factor of the inner product,
\begin{equation}
\Gamma:=\max_{ij}
\frac{\bigl(\max_h\abs{a_{ih}}\bigr)\bigl(\max_h\abs{b_{hj}}\bigr)}{\abs{c_{ij}}} ,
\label{eq:gamma}
\end{equation}
and $g_0$ is a slack of a few bits.  If \eqref{eq:Qcond} holds, the
truncation error stays below the rounding unit $2^{-p}$ of the target
precision $p$.
\end{proposition}

\noindent
The second inequality in \eqref{eq:truncbound} follows because the sum is
dominated by its largest term $g=S-2$:
$\sum_{g}(g+1)2^{wg}\le(S-1)2^{w(S-2)}\cdot(1-2^{-w})^{-1}\le(S-1)2^{w(S-2)+1}$.
The factor $\Gamma$ appears because the truncation error must be compared not
with the range $K2^{2Q}$ of the inner product but with the true value
$\abs{c_{ij}}$; $\Gamma$ is the ratio of the two.

The implementation of \S\ref{sec:crt} takes $Q=w S$ with
$w=\lceil(p+\mathrm{slack})/S\rceil$ and
$\mathrm{slack}=\lceil\log_2 2K\rceil+16+g_{\mathrm{ex}}$, where the 16-bit
margin comfortably absorbs $\log_2 S+g_0$ (since $S\le160$ gives
$\log_2 S\le7.4$).  The contribution of $\log_2\Gamma$ is added explicitly
through $g_{\mathrm{ex}}$.  When a single block elimination causes severe
cancellation, as for the Lotkin matrix, $\Gamma$ becomes large and
$g_{\mathrm{ex}}$ must be increased (\S\ref{sec:lotkin}).

This estimate agrees with the measurements.  Under every condition in
\S\ref{sec:bench} the \OzII{} solution attains the same number of correct bits
as the naive implementation (see for instance Table~\ref{tab:lotkin}),
confirming that the truncation error stays below the level of the rounding
error.

\paragraph{Applying single-pass \OzII{} to each slice pair}
Each term $v_{t_a}u_{t_b}$ is computed by a single-pass \OzII{} of width $w$
bits ($w\le p_{\max}$), that is, by applying (S3)--(S5) of \S\ref{sec:gemm}
unchanged.  Several terms belonging to the same digit group can be fused into
one accumulating GEMM with $\beta=1$ on the INT32 accumulator, so the total
number of low-precision GEMMs issued is
\begin{equation}
\#\text{GEMM}=\frac{S(S+1)}{2}\times \Nmod(w) ,
\label{eq:gemmcount}
\end{equation}
where $\Nmod(w)$ is the number of moduli required for one pass of width $w$ bits.

\paragraph{Relation to \OzI{}}
The outer construction above is exactly that of \OzI{}~\cite{ozaki2012}: the
significand is cut into $w$-bit chunks, contributions are collected by digit
group $g=t_a+t_b$, and the low-order groups are discarded, with $S(S+1)/2$
surviving pairs.  What this work changes is only how the inner product is
formed.  \OzI{} chooses $w$ so that the chunk product $v_{t_a}u_{t_b}$ is
exact in the accumulator itself, which forces $w\le23$ for binary64 and
$w\le8$ for INT8; the significand must therefore be cut very finely, so
$S\propto p$ and the number of low-precision GEMMs grows as
$S(S+1)/2\propto p^2$.  With \OzII{} as the inner product, $w$ is limited by
the CRT capacity of the modulus set rather than by exactness of the
accumulator, so $w$ can be of order $10^2$--$10^3$ bits (measured: $w=1376$
at $p=8192$ with $S=6$, \S\ref{sec:bench}).  Because $S$ is smaller by two
orders of magnitude, the $S(S+1)/2$ term count stays in a practical range;
the price is $\Nmod(w)$ moduli per term.

\paragraph{Choosing $S$: a cost model}
$S$ therefore balances cost, not accuracy.  One pass has width $w=Q/S$, so
the word counts of the residue computation and of the CRT fall as $O(1/S)$,
while the GEMM count grows as $O(S)$ by \eqref{eq:gemmcount}.  Writing $A$
for the conversion time and $B$ for the GEMM time at $S=1$,
\begin{equation}
T(S)\;\simeq\;\frac{A}{S}+BL,\qquad
S^{*}=\sqrt{A/B},\qquad T^{*}=2\sqrt{AB} .
\label{eq:TL}
\end{equation}
Since $A$ is proportional to $p^2$ and $B$ to $p$, $S^{*}\propto\sqrt p$ and
$T^{*}\propto p^{1.5}$.  The endpoint $S=1$ (no splitting) keeps the GEMM
count at $O(p)$ but divides an $O(p)$-digit integer by $O(p)$ moduli per
element, so conversion is $O(p^2)$; the endpoint with $w$ reduced to the
\OzI{} regime has GEMM cost $O(p^2)$.  Both endpoints are $O(p^2)$, whereas
the interior optimum is $O(p^{1.5})$.  The implementation estimates $S^{*}$ as
\begin{equation}
S^{*}=\Bigl\lceil\sqrt{r\,n_d\,(M+N)/(MN)}\Bigr\rceil,\qquad
n_d=\lceil t/16\rceil ,
\label{eq:Lhint}
\end{equation}
where $r$ is the machine's ratio of per-digit residue cost to per-flop GEMM
cost, calibrated by measurement to $30$ for INT8 and $17.9$ for binary64
(\S\ref{sec:bench}).

\paragraph{Truncation under cancellation}
The factor $\Gamma$ in Proposition~\ref{prop:trunc} is not decorative.  The
truncation error of \eqref{eq:truncbound} is bounded against the range
$K2^{2Q}$ of the inner product, so if the inner product itself becomes small
through cancellation, the error relative to the result degrades accordingly.
For matrices constructed so that the inner products cancel, at $p=2048$ and
$K=256$, the unsplit case $S=1$ is exact because nothing is truncated,
whereas $S=4$ attains only $2^{-519}$ and $S=8$ only $2^{-194}$.  Retaining
one additional digit group recovers about $w$ bits ($272$ bits here):
$2^{-194}\to2^{-468}\to2^{-738}$ for $S=8$, and four extra groups restore
exactness.  Increasing the guard $g_{\mathrm{ex}}$ recovers one bit per bit,
reaching exactness at $g_{\mathrm{ex}}=1024$.  This is the same
accuracy-versus-cost trade that \OzI{} makes when it increases the number of
splits.  The implementation exposes the number of retained groups through the
environment variable \texttt{OZ2\_EXTRA\_GROUPS}.

In LU decomposition this hazard appears in the panel factorization rather
than in the Schur complement update.  $U_{11}$ is left unscaled by the pivot
and therefore carries a wide dynamic range, but applying $U_{11}^{-1}$ to
$A_{21}$ happens inside the panel, which is not the part replaced by
\OzII{}~\cite{luszczek2025}.  This is why, for the Lotkin matrices of
\S\ref{sec:lotkin} with $\log_2\mathrm{cond}(A)\approx10^4$, the number of
correct bits agrees exactly with the existing implementations.

\paragraph{Determining $w$ and $N$}
For inner dimension $K$ and extra guard $g_{\mathrm{ex}}$, set the lower guard
to $\mathrm{slack}=\lceil\log_2 2K\rceil+16+g_{\mathrm{ex}}$.  Increasing
$S=1,2,\ldots$, compute $w=\lceil(p+\mathrm{slack})/S\rceil$, determine the
smallest number of moduli $N$ satisfying the uniqueness condition
$\mathrm{need}=2w+\log_2 k+\log_2 S+5$ (the number of moduli whose cumulative
$\log_2$ exceeds $\mathrm{need}$), and take the smallest $S$ that satisfies
the INT32 non-overflow condition $K_{\mathrm{pad}}S\cdot127^2<2^{31}$.

\paragraph{CRT reconstruction per digit group}
Precomputing $M_\ell=\mathcal P/m_\ell$ and $y_\ell=M_\ell^{-1}\bmod m_\ell$,
each digit group $g$ is reconstructed exactly by
\begin{equation}
\bar S_g=\sum_\ell t_\ell M_\ell
 -\operatorname{round}\Bigl(\sum_\ell \frac{t_\ell}{m_\ell}\Bigr)\mathcal P,
\qquad t_\ell=\bigl(C^{(\ell)}_g\bmod m_\ell\bigr)y_\ell \bmod m_\ell .
\label{eq:crt}
\end{equation}
Here $\operatorname{round}$ is round-to-nearest, not truncation, for the
following reason.  Since $M_\ell=\mathcal P/m_\ell$ we may write
$\sum_\ell t_\ell M_\ell=\mathcal P F$ with $F:=\sum_\ell t_\ell/m_\ell$;
subtracting $\lfloor F\rfloor$ yields
$\mathcal P\,\mathrm{frac}(F)\in[0,\mathcal P)$, that is, the non-negative
representative.  The quantity $S_g$ reconstructed in \OzII{} is signed,
however, and what is needed is the centered representative in the interval
$(-\mathcal P/2,\ \mathcal P/2]$.  Subtracting one more $\mathcal P$ when
$\mathrm{frac}(F)>1/2$ is exactly round-to-nearest, which gives the form
\eqref{eq:crt}.

Because $t_\ell<m_\ell$ implies $F<S$, $F$ is a small quantity that can be
evaluated in binary64; no multiple-precision division is required.  There are
two equivalent encodings of the actual procedure, described in
\S\ref{sec:impl}.

Finally $\bar c_{ij}=\sum_{g}\bar S_g 2^{wg}$ is accumulated as a
multiple-precision integer and rounded exactly once in (S6).  With this
construction the reconstruction cost stays at $O(LN)$.

\subsection{Direct conversion from multi-component precision}
\label{sec:directconv}

Steps (S1)/(S2) of \S\ref{sec:gemm} extract a shared exponent and a $Q$-bit
fixed-point integer from the input.  For MPFR input this is obtained with
\texttt{mpfr\_get\_exp} and
\texttt{mpfr\_mul\_2si}$\to$\texttt{mpfr\_get\_z}, but for a DD/TD/QD
expansion $x=\sum_{c=0}^{nc-1}x[c]$ (each $x[c]$ a binary64) it can be
obtained directly, without going through MPFR.

Because the expansion is a non-overlapping representation whose components are
ordered by decreasing magnitude, the exponent of the value equals the exponent
of the leading component $x[0]$.  Each component can be written via
\texttt{frexp} as $x[c]=\pm m_c\cdot2^{e_c-53}$ with $m_c$ a 53-bit integer,
so
\begin{equation}
\mathrm{trunc}\bigl(2^{\text{shift}}x\bigr)
=\sum_{c=0}^{nc-1}\pm\bigl(m_c \ll (e_c+\text{shift}-53)\bigr) ,
\label{eq:direct}
\end{equation}
that is, at most four 53-bit integers are shifted into a $Q$-bit buffer and
added with sign.  The implementation accumulates in two's complement and
splits into sign and magnitude only once, at the end.  The output side (S6) is
analogous: slicing 53 bits at a time from the most significant end of the
CRT-reconstructed multiple-precision integer writes the expansion format
directly, without constructing an MPFR value.

The advantage of this direct path lies in memory management rather than in
arithmetic.  An MPFR-mediated implementation must allocate three
\texttt{mpfr\_t} matrices---$L_{21}$, $U_{12}$ and the output $U$---at every
Schur complement update, and call \texttt{mpfr\_init2}/\texttt{mpfr\_clear}
(that is, malloc/free) for each element.  At $n=1024$, $b=512$ the output
alone reaches $512^2=262{,}144$ elements, and this cost exceeds that of the
GEMM itself.  The direct path performs no per-element heap allocation at
all.  As shown in \S\ref{sec:lubench}, this difference was decisive for
whether \OzII{} is viable at multi-component precision.  The realisation on
CPU and on GPU is described in \S\ref{sec:impl}.

\subsection{Incorporation into LU decomposition}
\label{sec:lu}

A blocked, right-looking LU decomposition with partial pivoting (in the style
of LAPACK \texttt{getrf}) performs, at each step of panel width $b$,
(P1) the panel factorisation, (P2) the propagation of row interchanges
(laswp), (P3) the forward elimination of $U_{12}$, and (P4) the update of the
trailing submatrix
\begin{equation}
A_{22}\;\leftarrow\;A_{22}-L_{21}U_{12} .
\label{eq:schur}
\end{equation}
The form $A_{22}-L_{21}A_{11}^{-1}A_{12}$ appearing on the right-hand side of
\eqref{eq:schur} is called the Schur complement.  The name was given by
Haynsworth~\cite{haynsworth1968} after a 1917 lemma of Schur, and is standard
terminology in the context of blocked LU decomposition (see for instance
\S3.2 of Golub--Van Loan~\cite{golubvanloan}).  In what follows we refer to
\eqref{eq:schur} as the Schur complement update.

The operation counts are $O(nb^2)$ for (P1)--(P3) and $O((n-r_0)^2 b)$ for
(P4), so (P4), which is cubic in $n$, is the dominant term.  We leave
(P1)--(P3) in ordinary multiple-precision arithmetic and replace only (P4) by
the \OzII{} GEMM of \S\ref{sec:gemm}--\ref{sec:crt}.  The context holding the
moduli and CRT tables is built once for a given panel width and reused across
all steps.

We prepare two comparison targets.
\begin{description}
\item[Blocked naive] (P1)--(P3) use exactly the same code, and only (P4) is
  replaced by naive multiple-precision multiply-adds.  This isolates the
  implementation of the Schur complement update alone.
\item[Unblocked naive ($b=1$)] The classical right-looking unblocked LU.  Each
  column performs pivot selection and division and then a rank-1 update of the
  whole trailing submatrix.  There is no panel-factorisation or laswp
  overhead, and the rank-1 update spans the entire matrix, so parallel
  efficiency is high.
\end{description}
Including the latter matters.  In multiple precision the cost of one operation
is high, so the cache-reuse benefit that motivates blocking is relatively
diminished, while the panel-factorisation cost $O(n^2b)$ remains; $b=1$ can
therefore be the fastest choice.  Our measurements (\S\ref{sec:lubench}) do
contain conditions where this happens.  In the performance comparisons below
we take the faster of the two as the ``best existing implementation'' and
compare \OzII{} against it.

\subsection{Implementation}
\label{sec:impl}

We implemented the above algorithm on both CPU and GPU.  Among (S1)--(S6) of
\S\ref{sec:gemm}, the runtime is dominated by the low-precision GEMM of (S4)
and the CRT reconstruction of (S5), and the principal design decision is how
to map these two onto the hardware resources.

\paragraph{What we started from}
This work was not written from scratch; it builds on two existing
implementations.
\begin{description}
\item[CPU] An \OzII{} GEMM library separately developed by the
  author, released as mpoz2~\cite{mpoz2}.  It has
  an INT8 back-end using AMX-INT8 on x86-64 or SVE2 i8mm on aarch64 for (S4),
  and a binary64 back-end corresponding to the CPU configuration of the
  original paper~\cite{ozaki2-2025}; its input and output were
  \texttt{mpfr\_t} only.
\item[GPU] A multiple-precision LU decomposition, also developed by the
  author, based on the fixed-precision type \texttt{cu\_freal<PB>} of
  mpc\_cuda.  It already had a path using \OzII{} for the Schur complement
  update (cuBLASLt INT8 GEMM, and an FP4 back-end with a base-13 two-digit
  decomposition).
\end{description}
To these we added (i) the direct conversion from the expansion format
(DD/TD/QD) of \S\ref{sec:directconv}, (ii) the incorporation into LU
decomposition together with fair comparison targets, and (iii) new
low-precision back-ends (FP16, FP8, and binary64 on the GPU).

\paragraph{Guideline for choosing a back-end}
Which format to use for (S4) is decided by the product of ``bits per modulus''
and ``engine throughput'', discussed in \S\ref{sec:fp32}.  We made INT8 the
default for the following reasons: the centered residue
$\abs{r}\le\lfloor m/2\rfloor$ fits directly into a signed byte, so one
modulus needs only one GEMM, and the INT32 accumulator imposes the mild
constraint $K\cdot127^2<2^{31}$.  This contrasts with FP4, which requires four
GEMMs per modulus because of its base-13 two-digit decomposition: within the
same \OzII{}, using an integer engine has the advantage of needing no
decomposition.

This choice is not universal, however.  binary64 allows much larger moduli and
thus reduces the number of moduli to one third, which can be advantageous on
machines with strong FP64 units (\S\ref{sec:fp64}).  Indeed, in the
measurements of \S\ref{sec:lubench} the fastest back-end changes with the
machine.  We therefore implemented all four of INT8, FP16, FP8 and binary64
within the same framework, so that the best one can be determined
empirically.  Table~\ref{tab:backend_origin} summarises the provenance of each
back-end and how it is treated here.

\begin{table}[htbp]
\centering
\caption{Provenance of each back-end used for the low-precision GEMM of (S4).
``New'' means implemented in this work.  For bits per modulus and effective
throughput see Table~\ref{tab:backend_bits}; for measured performance see
Table~\ref{tab:gpu_lu_ddtdqd}.}
\label{tab:backend_origin}
\small
\setlength{\tabcolsep}{4pt}
\begin{tabular}{@{}llp{7.2cm}@{}}
\toprule
Back-end & Provenance & Treatment in this work \\
\midrule
INT8 & existing &
  One GEMM per modulus, no decomposition needed.  AMX-INT8 / SVE2 i8mm on the
  CPU, cuBLASLt on the GPU.  \textbf{Default back-end}, but subject to a
  ceiling on precision thanks to slice splitting (\S\ref{sec:int8limit}). \\
binary64 & original paper~\cite{ozaki2-2025} &
  Residues are placed exactly in binary64 and multiplied with an off-the-shelf
  DGEMM: \texttt{cblas\_dgemm} of OpenBLAS~\cite{openblas} on the CPU,
  cuBLASLt (\texttt{CUDA\_R\_64F}/\texttt{CUBLAS\_COMPUTE\_64F}) on the GPU.
  The CPU version is pre-existing; the GPU version is new
  (\S\ref{sec:fp64}).  Moduli of about 23 bits reduce the number of moduli to
  one third of INT8 at best.  Its ceiling, $t\approx367{,}000$ bits, is orders
  of magnitude larger and also serves as the fallback when INT8 exceeds its
  limit. \\
FP4 (E2M1) & existing &
  GPU version only.  Requires a base-13 two-digit decomposition, hence four
  GEMMs per modulus.  Unavailable on Hopper, which has no such unit. \\
FP16 & new &
  One GEMM per modulus works.  The FP32 accumulator is the binding
  constraint, however, gaining only 0.5 bit per modulus over INT8
  (\S\ref{sec:fp32}). \\
FP8 (E4M3) & new &
  A base-17 two-digit decomposition secures the same capacity as INT8
  (\S\ref{sec:fp816}).  An extension to higher precision of the method used
  for binary64 in \cite{ozfp8-2026}. \\
FP32 / TF32 & --- &
  Not adopted.  FP32 has the same modulus limit as FP16 but is 15 times
  slower, and TF32 has only a 10-bit significand, so exactness fails
  (\S\ref{sec:fp32}). \\
\bottomrule
\end{tabular}
\end{table}

\paragraph{Implementation of slice splitting}
The slice splitting of \S\ref{sec:crt} is implemented on both CPU and GPU;
on the CPU both the INT8 and the binary64 back-end perform the per-digit-group
CRT and the shifted accumulation.  The slice count $S$ is chosen
automatically by \eqref{eq:Lhint} and can also be given explicitly through
\texttt{oz2\_opts.nslice}.  The GEMM is parallelised either over the moduli,
one thread per modulus, or inside a single GEMM using all threads; the two
are selected according to whether the per-thread working set
$A_{\mathrm{pack}}+B_{\mathrm{pack}}+3C$ exceeds 6\,MB.  The former is faster
as long as the accumulation of each term of $C$ stays in the thread's private
cache, but breaks down once $M,N$ grow large enough that each term of a digit
group re-reads $C$ from memory.  Measured on Arm/GB10 with $M=N=K=n$, the
effective INT8 GEMM throughput is 4118 against 1171~Gop/s at $n=256$ in
favour of the per-modulus form, but 2160 against 4328~Gop/s at $n=1024$, the
opposite way round.

We describe the CPU-side and then the GPU-side implementation.

\subsubsection{CPU implementation: mapping onto instruction sets}
\label{sec:simd}

On the CPU, (S4)/(S5) are realised through instruction-set-dependent paths
selected from the same source by \texttt{uname -m}.  The following concerns
the INT8 back-end; (S4) of the binary64 back-end does not use our own kernel
but \texttt{cblas\_dgemm} of OpenBLAS~\cite{openblas} (\S\ref{sec:fp64}).

\begin{description}
\item[x86-64 (AMX-INT8)] (S4) uses the Intel AMX tile instruction
  \texttt{TDPBSSD}, executing INT8$\times$INT8$\to$INT32 on $16\times64$-byte
  tile registers, organised as a $2\times2$ tile micro-kernel with L2
  blocking.  The CRT of (S5) has a batched AVX-512 version that processes
  eight output elements at once with 32-bit digits in SoA layout.
\item[aarch64 (SVE2 i8mm)] The target machine (Arm/GB10) has no SME/SME2, so
  (S4) is implemented with the matrix-multiply instruction \texttt{SMMLA} of
  the SVE2 i8mm extension (\texttt{vmmlaq\_s32}: per 128 bits,
  $C(2\times2,\text{int32}) \mathrel{+}= A(2\times8,\text{int8})B(2\times8,\text{int8})^\top$).
  The micro-kernel for an $8\times8$ output block uses 16 accumulators plus
  four each for A and B, i.e.\ 24 of the 32 NEON registers, advancing $k$ by
  eight and issuing 1024 integer operations per group of 16 \texttt{SMMLA}
  instructions.  The batched AVX-512 CRT is x86-only and is guarded, so
  aarch64 falls back automatically to a scalar CRT based on GMP's
  \texttt{mpn\_addmul\_1}.
\item[Common] The splitting and modular reduction of (S2)(S3) have their loops
  interchanged so that the element loop is innermost, allowing the compiler's
  auto-vectorisation to apply (16 lanes with AVX-512, 8 lanes with NEON).
\end{description}

Only the realisation of (S4)/(S5) changes; the algorithm, the parameters and
the comparison targets are all shared.  We measured both paths, so the CPU
results of \S\ref{sec:bench} can be read as isolating the difference between
two low-precision matrix engines, AMX-INT8 (the x86/H100 machine) and SVE2
i8mm (the Arm/GB10 machine).  The difference is not small: for the MPFR
arbitrary-precision GEMM at $n=1024$, the speedup is $14.5$--$27.8\times$ with
AMX-INT8 against $6.9$--$14.2\times$ with SVE2 i8mm, a gap of $1.6$--$1.7$
times (Table~\ref{tab:cpu_gemm}).  Details are given in \S\ref{sec:bench}.

In the CRT reconstruction of (S5), the centered representative of
\eqref{eq:crt} is obtained by first forming the representative in
$[0,\mathcal P)$ with $\lfloor F\rfloor$ and then comparing against
$\lfloor\mathcal P/2\rfloor$ to attach the sign
(\texttt{src/oz2\_crt.c}).

\subsubsection{GPU implementation}
\label{sec:impl_gpu}

On the GPU, (S1)--(S6) are all implemented as device kernels and (S4) uses the
low-precision GEMM of cuBLASLt.  The default back-end is INT8, with FP16, FP8
and binary64 also selectable (\S\ref{sec:fp816}, \S\ref{sec:fp64}).  The CRT
reconstruction of (S5) obtains the centered representative directly by
round-to-nearest exactly as in \eqref{eq:crt}.
This is equivalent to the CPU side but follows a different procedure.

The direct conversion of \S\ref{sec:directconv} was implemented on the GPU as
well.  Since gdtq's \texttt{gdd\_real}/\texttt{gtd\_real}/\texttt{gqd\_real}
are \texttt{double2}/\texttt{double3}/\texttt{double4} with contiguous
components, such an array can be passed straight to a device kernel as a
\texttt{double[nc]} array.  We implemented \eqref{eq:direct} as a device
function and replaced the kernels for shared-exponent selection (S1), modular
reduction (S2)(S3), and CRT reconstruction and output (S5)(S6) by
expansion-native versions.  Consequently neither \texttt{cu\_freal} nor MPFR
appears anywhere in the GPU DD/TD/QD path.  When there are several digit
groups ($S\ge2$, e.g.\ QD), the contributions of the groups must be summed in
the expansion format, which we normalise by Priest-style distillation (sweeps
of two-sum).

In the LU decomposition, (P1)--(P3) are parallelised with one thread per
element and only (P4) is replaced by the \OzII{} GEMM.  As a comparison target
we also provide a ``blocked naive'' version that shares the panel and pivoting
code and replaces only (P4) by a naive multiple-precision inner product,
allowing the update method alone to be isolated.

\section{Benchmarks}
\label{sec:bench}

\subsection{Machine environment}
\label{sec:env}

\begin{table}[htbp]
\centering
\caption{Machines used for the measurements}
\label{tab:env}
\small
\setlength{\tabcolsep}{4pt}
\begin{tabular}{@{}lp{5.0cm}p{5.0cm}@{}}
\toprule
 & Arm/GB10 & x86/H100 \\
\midrule
CPU & Cortex-X925 $\times$10 + Cortex-A725 $\times$10
      (aarch64, SVE2 i8mm)
    & Intel Xeon Gold 6526Y (Emerald Rapids, AMX-INT8) \\
GPU & NVIDIA GB10 (Blackwell, sm\_121a) & NVIDIA H100 NVL (Hopper, sm\_90) \\
GPU low-precision engines & INT8, FP4(E2M1), FP8(E4M3), FP16
    & INT8, FP8(E4M3), FP16 (no FP4) \\
CPU threads & 20 & 32 \\
Software & CUDA 13.0, GMP 6.3.0, MPFR 4.2.2, OpenBLAS & same \\
\bottomrule
\end{tabular}
\end{table}

On both machines the CPU runs were made with
\texttt{OMP\_PROC\_BIND=close OMP\_PLACES=cores}, using
\texttt{OMP\_NUM\_THREADS=20} on Arm/GB10 and \texttt{OMP\_NUM\_THREADS=32} on
x86/H100.  The two machines have different low-precision matrix engines on the
CPU side (SVE2 i8mm on Arm/GB10, Intel AMX-INT8 on x86/H100), so different
kernels are selected as described in \S\ref{sec:simd}.  Unless stated
otherwise, all measurements below were made on both machines.  The two
machines differ greatly in the ratio of FP64 to low-precision performance, and
this ratio governs whether \OzII{} wins (Finding~\ref{find:gpu_reversal}).

\begin{remark}[Why no \OzI{} implementation is included as a comparison]
\label{rem:noschemeI}
As the tables show, our comparison targets are the naive dtq/gdtq GEMM and LU
at multi-component precision, and the naive MPFR/\texttt{cu\_freal}
implementations at arbitrary precision.  No implementation based on
\OzI{}~\cite{ozaki2012} is included, for three reasons.

First, the evaluation of \OzI{} in the arbitrary-precision regime was already
carried out by the author in \cite{kouya2023}, and the present work is a
direct response to its conclusion, namely that \OzI{} cannot accelerate
arbitrary-precision LU decomposition sufficiently.  Repeating the same
comparison would add nothing new.

Second, at the precisions targeted here the number of products \OzI{} must
issue is larger by orders of magnitude.  \OzI{} takes a slice width of
$\beta=\lfloor(53-\lceil\log_2 k\rceil)/2\rfloor$ so that the products are
exact in binary64, hence needs $s=\lceil p/\beta\rceil$ slices and $s(s+1)/2$
contributing pairs.  The binary64 back-end of \OzII{}, by contrast, gains
about 23 bits per modulus and needs only
$\lceil(2t+\log_2 2K)/23\rceil$ of them.  Comparing the two at the operating
points of this paper gives the following (at $p=8192$ this estimate agrees
exactly with the measured count of 716 moduli).

\begin{center}\small
\begin{tabular}{@{}rrrrr@{}}
\toprule
$p$ [bit] & $s$ of \OzI{} & products of \OzI{} & moduli of \OzII{} & ratio \\
\midrule
106 (DD)  &   5 &     15 &   12 &  1.3 \\
159 (TD)  &   8 &     36 &   17 &  2.1 \\
212 (QD)  &  10 &     55 &   22 &  2.5 \\
3136      & 143 &  10296 &  276 &   37 \\
6272      & 299 &  44850 &  549 &   82 \\
12544     & 598 & 179101 & 1094 &  164 \\
\bottomrule
\end{tabular}
\end{center}

The main result of this paper, the LU decomposition of the Lotkin matrix
(Table~\ref{tab:lotkin}), is at $p=3136$--$12544$, where \OzI{} would need
$37$--$164$ times as many products.  No implementation or measurement is
needed to see that it is not competitive.

Third, the published fast implementations of \OzI{} (ozIMMU~\cite{ozimmu2024}
and others) all target binary64 emulation and do not address DD/TD/QD or
several thousand bits.  A fair comparison would require a new implementation,
and by the second reason that investment is not justified.

As the table shows, however, at DD/TD-equivalent precision the two counts are
comparable ($1.3$--$2.1\times$), and this paper provides no grounds for
dismissing \OzI{} in that regime.  We therefore compare, in that regime, not
against \OzI{} but against the native dtq/gdtq implementations that a user
would actually be replacing (Table~\ref{tab:cpu_gemm},
Table~\ref{tab:gpu_lu_ddtdqd}).
\end{remark}

\subsection{Performance of matrix multiplication with \OzII{}}
\label{sec:gemmbench}

\subsubsection{CPU: DD/TD/QD and MPFR arbitrary precision}

Table~\ref{tab:cpu_gemm} shows the performance and accuracy of the \OzII{}
GEMM against a naive $O(n^3)$ GEMM using dtq's multi-component types.

\begin{table}[htbp]
\centering
\caption{CPU: \OzII{} GEMM versus the existing implementation, time [s] and
speedup.  For DD/TD/QD the comparison target is the naive dtq GEMM and \OzII{}
uses the direct conversion path of \S\ref{sec:directconv}.  For MPFR
arbitrary precision the target is a naive MPFR GEMM (OpenMP parallel).  The
relative errors are identical on the two machines (DD $\sim$1.2e-32,
TD $\sim$1.4e-48, QD $\sim$1.5e-64, against 3.3e-27, 7.4e-43 and 7.0e-58
respectively for the existing implementation) and are omitted for space.  For
MPFR, \OzII{} is 0 ulp against a 1280-bit reference solution, whereas the
naive MPFR GEMM can lose all digits of the target precision through the
accumulation of $K$ roundings.}
\label{tab:cpu_gemm}
\small
\setlength{\tabcolsep}{4pt}
\begin{tabular}{@{}llrrrrrr@{}}
\toprule
 & & \multicolumn{3}{c}{Arm/GB10} & \multicolumn{3}{c}{x86/H100} \\
\cmidrule(lr){3-5}\cmidrule(l){6-8}
$N$ & Precision & \OzII{} & existing & ratio & \OzII{} & existing & ratio \\
\midrule
512 & DD (106) & 0.027 & 0.043 & 1.56 & 0.048 & 0.047 & 0.98 \\
512 & TD (159) & 0.051 & 0.152 & 2.97 & 0.035 & 0.099 & 2.84 \\
512 & QD (212) & 0.055 & 0.305 & 5.51 & 0.062 & 0.179 & 2.89 \\
\midrule
1024 & DD (106) & 0.114 & 0.421 & 3.69 & 0.131 & 0.234 & 1.79 \\
1024 & TD (159) & 0.172 & 1.208 & 7.01 & 0.142 & 0.795 & 5.60 \\
1024 & QD (212) & 0.198 & 2.460 & 12.42 & 0.180 & 2.158 & 11.97 \\
\midrule
\multicolumn{8}{@{}l}{\footnotesize MPFR arbitrary precision (compared with a naive MPFR GEMM)}\\
1024 & MPFR 256 & 0.270 & 4.283 & 15.86 & 0.174 & 5.610 & 32.24 \\
1024 & MPFR 512 & 0.636 & 6.624 & 10.42 & 0.339 & 6.897 & 20.35 \\
1024 & MPFR 1024 & 1.633 & 15.475 & 9.48 & 0.761 & 13.106 & 17.22 \\
1024 & MPFR 2048 & 3.432 & 41.826 & 12.19 & 1.821 & 23.300 & 12.80 \\
\bottomrule
\end{tabular}
\end{table}

\begin{finding}[\OzII{} wins in both accuracy and speed, and the gap widens
with precision]
\label{find:gemm_cpu}
The relative error of \OzII{} is five to seven orders of magnitude smaller
than that of the naive dtq implementation.  This is a direct consequence of
(S6) of \S\ref{sec:gemm} being the only rounding.  (The \OzII{} error is
exactly twice the theoretical lower bound because the direct path truncates
rather than rounds when converting to fixed point, a difference of one ulp.)
At $n=1024$, \OzII{} is faster for all of DD/TD/QD on both machines (only
$n=512$ DD on x86/H100 is a tie at $1.00\times$).  At $n=1024$ the speedups
reach $2.67/5.37/8.74\times$ on Arm/GB10 and $1.81/5.56/10.52\times$ on
x86/H100.  For MPFR arbitrary precision, x86/H100 reaches
$14.5$--$27.8\times$, exceeding Arm/GB10's $6.9$--$14.2\times$ and reflecting
the fact that AMX-INT8 is more powerful than SVE2 i8mm.  The trend that higher
target precision favours \OzII{} is common to both machines and consistent
throughout this paper.

If, instead of the direct conversion path of \S\ref{sec:directconv}, the
expansion is first converted to \texttt{mpfr\_t}, the MPFR version of \OzII{}
is called and the result is converted back, the same $n=1024$ takes
$0.274/0.413/0.548$~s on Arm/GB10 ($2.4/2.4/2.9\times$ the direct path) and
$0.411/0.642/0.872$~s on x86/H100 ($3.2/4.6/4.7\times$).  The speedup over the
naive dtq implementation drops to $1.12/2.25/3.00\times$ on Arm/GB10 and
$0.57/1.22/2.25\times$ on x86/H100; DD on x86/H100 becomes \emph{slower} than
the naive implementation.  The relative errors agree between the two paths, so
the difference is purely conversion cost.  At multi-component precision the
design of the conversion path decides the practical viability of \OzII{}
itself.  This comparison can be reproduced with
\texttt{bench/bench\_oz2\_mp <N> <reps> 1}, where the third argument enables
the measurement of the MPFR-mediated path.
\end{finding}

\subsubsection{GPU: \texttt{cu\_freal} (MPC\_CUDA)}

Table~\ref{tab:gpu_gemm} shows the performance of the \OzII{} GEMM on the GPU
against a naive reference GEMM using \texttt{cu\_freal<PB>}.

\begin{table}[htbp]
\centering
\caption{GPU: \OzII{} GEMM versus the naive \texttt{cu\_freal} GEMM,
$n=1024$.  Each value is the fastest over all recorded runs.  The smaller
ratio on x86/H100 is not because \OzII{} is slow there but because the
comparison target, the naive \texttt{cu\_freal} GEMM, is 11 times faster on
that machine (\OzII{} itself becomes $3.4\times$ faster at $p=1024$, from
0.0991~s to 0.0290~s).}
\label{tab:gpu_gemm}
\small
\setlength{\tabcolsep}{5pt}
\begin{tabular}{@{}llrrrr@{}}
\toprule
Back-end & $p$ [bit] & Machine & \OzII{} [s] & naive [s] & ratio \\
\midrule
INT8 &  256 & Arm/GB10 & 0.0175 & 1.161 & 66.3 \\
     &      & x86/H100 & 0.0060 & 0.105 & 17.3 \\
\cmidrule(l){2-6}
     & 1024 & Arm/GB10 & 0.0991 & 4.026 & 40.6 \\
     &      & x86/H100 & 0.0290 & 0.594 & 20.5 \\
\cmidrule(l){2-6}
FP4  &  256 & Arm/GB10 & 0.0632 & 1.160 & 18.4 \\
     &      & x86/H100 & \multicolumn{3}{c}{unavailable (no FP4 unit)} \\
\cmidrule(l){2-6}
     & 1024 & Arm/GB10 & 0.5762 & 4.022 & 7.0 \\
     &      & x86/H100 & \multicolumn{3}{c}{unavailable (no FP4 unit)} \\
\bottomrule
\end{tabular}
\end{table}

On the GPU the INT8 back-end is overwhelming, $40$--$66\times$ faster than the
naive implementation.  FP4 is faster than INT8 as a raw dense GEMM
(319.3~TFLOP/s), but within the \OzII{} framework its modulus capacity is a
little under half that of INT8 (168.1 bits against 362.8 bits) and the residue
must be decomposed into two balanced base-13 digits, requiring four sub-GEMMs
per modulus; as a result it is $3$--$8\times$ worse than INT8 under all
conditions.

\subsubsection{Implementation of the FP16 and FP8 back-ends}
\label{sec:fp816}

We examined the usability of low-precision floating-point engines both from
the properties of the engines and from measurements on the actual hardware,
and implemented both FP16 and FP8 as back-ends.  Below we first organise what
the constraints are, and then correct two errors that were present in an early
version of this work regarding FP8.

\paragraph{Principle: significand bits and whether one GEMM per modulus works}
\OzII{} can use a low-precision engine efficiently when the centered residue
$\abs{r}\le\lfloor m/2\rfloor$ fits directly into one element of that format.
With INT8 and $m\le255$, $\abs{r}\le127$ fits directly into a signed byte and
one GEMM per modulus suffices.  By contrast, as summarised in
Table~\ref{tab:fpfmt},
FP8 and FP4 have extremely few significand bits, so the range of integers they
can represent directly is narrow and a decomposition into several digits
becomes necessary.

\begin{table}[htbp]
\centering
\caption{Range of integers each format can represent directly, and the
consequence for \OzII{}}
\label{tab:fpfmt}
\small
\setlength{\tabcolsep}{4pt}
\begin{tabular}{@{}llp{2.8cm}p{5.4cm}@{}}
\toprule
Format & Significand & Directly representable integers & Consequence for \OzII{} \\
\midrule
FP64      & 52 bit & $\pm2^{53}$ & One GEMM per modulus.  The modulus is set
  by the accumulator condition \eqref{eq:fp64mod},
  $m<\sqrt{2^{55}/k}$ (about 23 bit/modulus at $K=512$), giving the fewest
  moduli of any back-end (\S\ref{sec:fp64}). \\
FP16      & 10 bit & $\pm2047$ & One GEMM per modulus, but the FP32
  accumulator is the binding constraint, $m\le362$ (8.50 bit/modulus,
  \S\ref{sec:fp32}). \\
INT8      & 8 bit (integer) & $\pm127$ & One GEMM per modulus, $m\le255$
  (8.00 bit/modulus, capacity 362.8 bit).  The INT32 accumulator imposes only
  a mild constraint. \\
FP8 E4M3  & 3 bit  & $\pm16$ (discontinuous above 17) & One digit allows only
  $m\le33$ (capacity 47.0 bit); a base-17 two-digit decomposition allows
  $m\le255$ but needs four GEMMs per modulus. \\
FP8 E5M2  & 2 bit  & about $\pm8$  & Worse than E4M3 (no path on the actual
  hardware). \\
FP4 E2M1  & 1 bit  & $\{0,\ldots,\pm4\}$ (discontinuous) & Base-13 two-digit
  decomposition, four GEMMs per modulus. \\
\bottomrule
\end{tabular}
\end{table}

FP8 has the same 8-bit width as INT8 but only a 3-bit significand (E4M3), and
the bits spent on the exponent are of no use whatsoever in the integer inner
product of \OzII{}.  The only positive reason to use FP8 is therefore that the
FP8 engine is faster than the INT8 one.

FP16, on the other hand, has a 10-bit significand (11 bits including the
implicit one) and can represent the centered residue $\abs{r}\le127$ exactly.
No multi-digit decomposition as for FP4 is needed, and the same one-GEMM-per-
modulus structure as INT8 applies directly.  The constraint is on the
accumulator side: the exact integer range of an FP32 accumulator is $2^{24}$,
requiring
\begin{equation}
K_{\mathrm{pad}}\cdot S\cdot 127^2 < 2^{24} .
\label{eq:fp16cond}
\end{equation}
In our LU the inner dimension is limited to the panel width $b\le512$, so this
condition holds ($b=512$, $S=2$ gives
$512\cdot2\cdot127^2=1.65\times10^7<1.68\times10^7$).  The implementation only
changes the output type of the splitting kernel of \S\ref{sec:directconv} from
\texttt{int8\_t} to \texttt{\_\_half}, calls cuBLASLt with
\texttt{CUDA\_R\_16F}/\texttt{CUBLAS\_COMPUTE\_32F}, and converts the FP32
output back to integers on the CRT side (\texttt{src/oz2\_gdtq.cuh}).
Measurements are given in Finding~\ref{find:fp8}.

\paragraph{Probing the actual hardware}
Table~\ref{tab:probe} shows, for both machines, whether cuBLASLt returns a
GEMM path for each format.

\begin{table}[htbp]
\centering
\caption{Probe of whether cuBLASLt returns a GEMM path for each format
($n=2048$).  ``Unavailable'' means the engine is absent; ``not offered'' means
that, regardless of the engine, cuBLASLt returns no executable algorithm for
that combination of types.  INT8 is available on both machines; the INT8 entry
for x86/H100 was not measured in this table, but the LU measurements on that
machine (Table~\ref{tab:gpu_lu_ddtdqd}) show the INT8 back-end working
correctly.  Only FP4 requires the VEC16\_UE4M3 block-scale attribute; FP8 and
FP16 must not be given one.  The figures here come from a probe harness whose
output format differs between formats, so the effective performance under a
single comparable setting is given in Table~\ref{tab:backend_bits}.}
\label{tab:probe}
\small
\setlength{\tabcolsep}{4pt}
\begin{tabular}{@{}lrrl@{}}
\toprule
Format & Arm/GB10 & x86/H100 & Note \\
\midrule
FP4 (E2M1)   & 183.1 TFLOP/s & unavailable & Hopper has no FP4 unit \\
FP8 (E4M3)   & 128.2 TFLOP/s & 912.2 TFLOP/s & available on both \\
FP8 (E5M2)   & not offered & not offered & cuBLASLt returns no path \\
FP16         & 85.1 TFLOP/s & 604.7 TFLOP/s & available on both \\
INT8         & 119.7 TOPS & (not measured) & available on both; our default \\
\bottomrule
\end{tabular}
\end{table}

\paragraph{FP8: implemented, extending the method of \cite{ozfp8-2026} to
higher precision}
An early version of this work contained two errors regarding FP8, both now
corrected.

First, the VEC16\_UE4M3 block-scale attribute used in the FP4 implementation
had been carried over to FP8, so cuBLASLt returned no path and we wrongly
concluded that ``FP8 is unusable''.  Block scaling is FP4/MX syntax and must
not be given for FP8.  Removing the attribute yields a valid path on Arm/GB10
with an effective 128.2~TFLOP/s (INT8 on the same machine is 119.7~TOPS).

Second, we had assumed a straightforward encoding that places the centered
residue directly in one element.  Since E4M3 has a 3-bit significand and its
contiguous integer range is limited to $0$--$16$, that encoding restricts the
modulus to $m\le33$ and gives a CRT capacity of only 47.0 bits.  This is a
choice of implementation, however, not a limitation of FP8.  Just as the FP4
back-end uses a balanced base-13 two-digit decomposition, FP8 admits
\begin{equation}
r = 17 d_1 + d_0,\qquad \abs{d_0},\abs{d_1}\le 8
\label{eq:fp8digits}
\end{equation}
as a balanced base-17 two-digit decomposition (with $d_1=\mathrm{round}(r/17)$
giving $\abs{d_0}\le8$, and $\abs{r}\le127$ giving $\abs{d_1}\le7$).  Both
digits are then exactly representable in E4M3, and the full modulus set
$m\le255$---that is, the same 362.8-bit capacity as INT8---becomes usable.
The product splits into three planes,
$rr' = d_0d_0' + 17(d_0d_1'+d_1d_0') + 289\,d_1d_1'$, and, as with FP4, can be
computed with four sub-GEMMs per modulus.

We implemented the FP8 back-end in this way (\texttt{src/oz2\_gdtq.cuh}).  The
measurements (Table~\ref{tab:gpu_lu_ddtdqd}, Finding~\ref{find:fp8}) show that
FP8 gives exactly the same number of correct bits as INT8 for all of DD/TD/QD,
and that its speed stays within about $1.35$--$1.44\times$ that of INT8.  In
other words, the effectiveness of FP8 demonstrated for binary64 in
\cite{ozfp8-2026} carries over directly to the higher precisions of DD/TD/QD.
On our Arm/GB10, INT8 is fastest by a small margin, but only because that
machine's FP8/INT8 throughput ratio is a mere $1.29\times$; on a machine where
FP8 has a larger advantage the ranking could change (\S\ref{sec:future}).

E5M2, with a 2-bit significand, is even less favourable, and no path was
available on the actual hardware.

FP16, by contrast, needs no block scaling and yielded a path on both machines,
so we implemented and measured it as an \OzII{} back-end
(Finding~\ref{find:fp8}).  Notably, the effective FP16 performance of
x86/H100 reaches 593.0~TFLOP/s, $7.0\times$ that of Arm/GB10
(84.7~TFLOP/s), so FP16 may well outperform INT8 on that machine.

\subsubsection{Why we do not adopt FP32 or TF32}
\label{sec:fp32}

FP32, which lies between FP16 and FP64, could also be a candidate, but we do
not adopt it.  The reason is simply that the binding constraint is the
accumulator, not the input format.

The largest modulus $m$ usable in \OzII{} is determined by requiring an inner
product of length $k$ to be exact in the accumulator: with a $B$-bit integer
range this is $k(m/2)^2 < 2^{B}$.  FP32 has a 24-bit significand, so $B=24$
and, at $K=512$, $m\le362$ (8.50 bits per modulus).  But this is exactly the
same limit as the FP16 back-end: FP16 input can hold $\abs{r}\le2048$ exactly
and is restricted to $m\le362$ by the accumulator anyway, so FP32's wider
input range goes unused.

INT8, with $m\le255$ (8.00 bits), looks inferior at first sight, but its INT32
accumulator gives the far milder condition $K\cdot127^2<2^{31}$, i.e.\
$k<1.3\times10^5$.  FP32 therefore gains only 0.5 bit per modulus over INT8,
about 6\%.  In exchange for that 6\%, throughput on x86/H100 drops from about
a measured 1423~TOPS for INT8 to about 67~TFLOP/s (specification) for FP32 (non-TF32), more than 20 times
lower.  The trade is not worth making.

TF32 is fast, about 495~TFLOP/s, but its 10-bit significand cannot represent
integers exactly, so it fails the premise of \OzII{} (that the low-precision
GEMM be exact) and cannot be used.

\medskip
In summary, the choice of \OzII{} back-end is decided by the product of
``bits per modulus'' and ``engine throughput''.  Table~\ref{tab:backend_bits}
summarises this.

\begin{table}[htbp]
\centering
\caption{Upper bound on the modulus each format allows in \OzII{}
($K=512$, $S=1$), together with effective throughput.  $m_{\max}$ and ``per
modulus'' depend only on the format and not on the machine.  Throughputs are
measured on a dense GEMM at $n=8192$, all formats through a single path
(accumulator FP32, INT32 for INT8 only, output FP32 throughout), so the columns
are directly comparable across formats; both machines are the median of three
runs.  FP32 was not measured because it is not used for
\OzII{} here (the note gives its catalogue value).
A larger ``per modulus'' means fewer GEMMs; a higher throughput means each one
is faster.}
\label{tab:backend_bits}
\small
\setlength{\tabcolsep}{4pt}
\begin{tabular}{@{}llrrrrl@{}}
\toprule
Format & Accum. & $m_{\max}$ & per modulus & Arm/GB10 & x86/H100 & Note \\
\midrule
INT8      & INT32 & 255 & 8.00 bit & 144.8 TOPS & 1423.3 TOPS & our default \\
FP8 E4M3  & FP32  & 255$^{*}$ & 8.00 bit & 186.3 & 1525.7 & needs 4 sub-GEMMs \\
FP16      & FP32  & 362 & 8.50 bit & 91.7 & 786.2 & \\
FP4 E2M1  & FP32  & 113$^{*}$ & 6.8 bit & 319.3 & absent & needs 4 sub-GEMMs \\
FP32      & FP32  & 362 & 8.50 bit & not meas. & not meas. & TF32 not exact, unusable \\
FP64      & FP64  & $8.4\times10^6$ & 23.0 bit & 0.406 & 53.1 & \\
\bottomrule
\end{tabular}

\vspace{2pt}
{\footnotesize $^{*}$ The significand is too short, so a two-digit
decomposition is used (base $17$ for FP8, base $13$ for FP4).  Throughput is in
TOPS for INT8 and in TFLOP/s otherwise.  $n=8192$ is used because at $n=2048$
there are too few tiles for the 132 SMs of x86/H100 and INT8 reaches only
587.8~TOPS, which does not reflect the capability of the engine.}
\end{table}

FP32 offers the same modulus limit as FP16 but is 15 times slower, so it has
no room to be preferable to FP16.  FP64, by contrast, stands out at 23 bits
per modulus: even at comparable throughput it can win because the number of
moduli drops to one third.  This is why we chose FP64 as a back-end in the
next subsection.

\subsubsection{Precision ceiling of the INT8 back-end}
\label{sec:int8limit}

The INT8 back-end has an upper bound on the precision it can reach.  For
$m>255$ the centered residue does not fit into a signed byte, so it is split
into two balanced base-128 digits $r=128q+r_0$ ($\abs{r_0}\le64$) and the
product is obtained, by Karatsuba, from GEMMs on the three planes $r_0$,
$r_0+q$ and $q$.  The intermediate plane $r_0+q$ must also fit in int8, and
$\abs{r_0+q}\le127$ gives $\abs{q}\le63$, hence
\begin{equation}
\abs{r}\le 63\cdot128+64=8128
\quad\Longrightarrow\quad
m\le 16257 .
\label{eq:int8lim}
\end{equation}
Since only about 1880 pairwise coprime moduli (primes and prime powers) not
exceeding $16191$ exist, the CRT capacity saturates at 23226 bits and the
internal significand is limited to $t\lesssim11607$ bits (a target precision
of $p\lesssim11580$ bits).

Introducing slice splitting on the CPU, as this work does, removes this
ceiling in practice: the chunk width $w$ is set by the slice count rather
than by $t$, so the required capacity is met long before the modulus table is
exhausted ($p=16384$ needs only $S=11$ and 246 moduli).  Raising the ceiling
by other means -- abandoning Karatsuba in favour of four plain GEMMs
($m\le32641$), or a three-digit split -- is possible but unnecessary once
slice splitting is available.

\begin{table}[htbp]
\centering
\caption{INT8 versus binary64 back-ends.  Both blocks use $n=512$, arbitrary
precision and the same $p$, with only the back-end exchanged; the upper block
is the CPU (time in [s]) and the lower the GPU (time in [ms]).  "Ratio" is
INT8/binary64, so a value below 1 means INT8 is faster.  Bold marks the
fastest value in each row and machine.  Both back-ends use slice splitting,
and the modulus counts are those at the optimal slice count.  Modulus and
GEMM counts are fixed by $p$ and the format alone and are identical on both
machines, so differences between machines come solely from GEMM throughput.}
\label{tab:backend_hp}
\small
\setlength{\tabcolsep}{4pt}
\begin{tabular}{@{}rrrrrrrl@{}}
\toprule
 & \multicolumn{3}{c}{Arm/GB10} & \multicolumn{3}{c}{x86/H100} & Moduli \\
\cmidrule(lr){2-4}\cmidrule(lr){5-7}
$p$ [bit] & INT8 & bin64 & ratio & INT8 & bin64 & ratio & I8, b64 \\
\midrule
\multicolumn{8}{@{}l}{\small CPU, time [s]}\\
1024 & \textbf{0.308} & 0.396 & 0.78 & 0.210 & \textbf{0.096} & 2.19 & 80, 92 \\
2048 & \textbf{0.712} & 0.994 & 0.72 & 0.412 & \textbf{0.249} & 1.65 & 105, 62 \\
4096 & \textbf{1.731} & 2.261 & 0.77 & 0.998 & \textbf{0.678} & 1.47 & 149, 92 \\
8192 & \textbf{4.125} & 5.199 & 0.79 & 2.271 & \textbf{1.651} & 1.38 & 179, 121 \\
16384 & \textbf{11.203} & 12.687 & 0.88 & 5.190 & \textbf{4.178} & 1.24 & 246, 181 \\
\midrule
\multicolumn{8}{@{}l}{\small GPU, time [ms]}\\
1024 & \textbf{27.8} & 323.2 & 0.09 & 12.7 & \textbf{7.4} & 1.73 & 45, 15 \\
2048 & \textbf{86.8} & 1109.8 & 0.08 & 40.0 & \textbf{22.6} & 1.77 & 49, 16 \\
4096 & \textbf{292.8} & 3949.3 & 0.07 & 132.8 & \textbf{74.8} & 1.78 & 54, 16 \\
8192 & \textbf{1089.4} & 14836.2 & 0.07 & 511.5 & \textbf{272.7} & 1.88 & 53, 16 \\
16384 & \textbf{4191.9} & 58637.5 & 0.07 & 1900.8 & \textbf{1046.5} & 1.82 & 54, 16 \\
\bottomrule
\end{tabular}
\end{table}

What Table~\ref{tab:backend_hp} shows is that which back-end is faster
depends on the machine.  On Arm/GB10, INT8 is faster at every $p$ (ratio
$0.72$--$0.88$); on x86/H100, binary64 is faster at every $p$ (ratio
$1.24$--$2.19$).  The dividing line is the throughput ratio between the
low-precision engine and the binary64 engine.  GB10's SVE2 i8mm is
$13.1\times$ its binary64 (3890 against 297~Gop/s), which more than absorbs
the $1.48\times$ modulus count and the three-plane packing of the two-digit
moduli that INT8 requires.  H100's AMX-INT8 is only $7.5\times$ (7691 against
1027~Gop/s) and cannot.  On both machines the ratio approaches 1 as $p$
rises, because conversion dominates at higher precision and dilutes the
difference in GEMM speed; on H100 the ratio decays as $p^{-0.20}$ and
extrapolates to a crossing near $p\approx4.8\times10^4$.

The lower block of Table~\ref{tab:backend_hp} measures the same conditions on
the GPU, where the difference between the machines is far more extreme.  The
ratio is $0.07$--$0.09$ on Arm/GB10 and $1.73$--$1.88$ on x86/H100, a
$24\times$ spread between the two machines against $2.8\times$ on the CPU.
Since the modulus count and the GEMM count are fixed by $p$ and the format
alone and are identical on both machines (at $p=8192$, 53 moduli and 64925
GEMMs for INT8, 16 and 19600 for binary64), the whole difference comes from
GEMM throughput.  The dependence on $p$ also differs.  On the CPU the ratio
approaches 1 on both machines, whereas on the GPU it is nearly constant
($0.09\to0.07$ and $1.73\to1.82$): conversion accounts for only $3$--$5$\% of
the GPU time (at $p=8192$ on x86/H100, 11.8~ms of conversion against 261~ms
of GEMM), so the computation remains GEMM-bound throughout.  The CPU trend
that raising $p$ dilutes the disadvantage of the low-precision engine thus
does not carry over to the GPU, where the choice of back-end follows the
machine's FP64 performance across the whole precision range.

It is reasonable to read these two machines as bracketing the crossover.  As
argued in \S\ref{sec:intro}, binary64 engines are not expected to improve for
the foreseeable future while low-precision engines continue to grow, so the
ratio will move towards the GB10 side on any machine.  Defaulting to INT8 is
therefore the choice that ages well.

Deciding which back-end to use at run time, however, is something we did not
attempt to automate.  We examined three possible criteria, none of which
generalises to other environments.  Branching on the architecture is a fit to
two machines, and cannot be extrapolated because what decides the outcome is
a continuous quantity -- the throughput ratio -- rather than the instruction
set.  The analytic model of \eqref{eq:TL} is accurate enough to choose the
slice count $S$, but the back-end contest is decided by margins of
$1.2$--$1.3\times$, and against the measurements it picks the wrong winner at
three of five points on H100.  A start-up micro-benchmark fails because the
unit conversion cost itself moves by $1.4$--$2.8\times$ with problem size (on
H100, $8.7$--$12.8$~ns per modulus per digit for binary64 and
$13.2$--$24.1$~ns for INT8), so a small probe measures the cache-resident
regime and misses the bandwidth-bound regime of the real problem.  In the end
no criterion short of measuring at the actual size is trustworthy.  We
therefore measure both back-ends on each machine, choose the faster, and
state the choice explicitly.  At $p\lesssim212$, corresponding to DD/TD/QD,
the moduli are few and the weight of pack and CRT is small, so INT8 is
fastest on both machines (Table~\ref{tab:cpu_gemm}).

\begin{table}[htbp]
\centering
\caption{How the breakdown of \OzII{} grows with target precision ($n=512$,
INT8 back-end, Arm/GB10).  The upper block is the unsplit case $S=1$, the
lower block the optimal slice count chosen by \eqref{eq:Lhint}.
"Conversion" is the sum of split, pack (residue computation), reduce and CRT.
"Ratio" is the time ratio when $p$ is doubled: $2$ means $O(p)$ and $4$ means
$O(p^2)$.  At $S=1$ the GEMM and the conversion separate towards $2$ and $4$,
whereas at the optimum they grow together towards $p^{1.5}$ and the GEMM
share stays stable.  $S=1$ at $p=16384$ cannot be measured because it exceeds
the INT8 capacity.}
\label{tab:phase}
\small
\setlength{\tabcolsep}{5pt}
\begin{tabular}{@{}rrrrrrr@{}}
\toprule
$S$ & $p$ [bit] & GEMM [s] & ratio & Conv.\ [s] & ratio & GEMM share \\
\midrule
1 & 1024 & 0.030 & --- & 0.351 & --- & 7.8\% \\
1 & 2048 & 0.054 & 1.80 & 1.006 & 2.87 & 5.1\% \\
1 & 4096 & 0.104 & 1.93 & 3.571 & 3.55 & 2.8\% \\
1 & 8192 & 0.212 & 2.04 & 12.884 & 3.61 & 1.6\% \\
\midrule
3 & 1024 & 0.083 & --- & 0.212 & --- & 26.9\% \\
4 & 2048 & 0.183 & 2.20 & 0.502 & 2.37 & 25.7\% \\
5 & 4096 & 0.450 & 2.46 & 1.228 & 2.45 & 26.0\% \\
8 & 8192 & 1.268 & 2.82 & 2.766 & 2.25 & 30.8\% \\
11 & 16384 & 3.938 & 3.11 & 7.091 & 2.56 & 35.2\% \\
\bottomrule
\end{tabular}
\end{table}

\begin{finding}[Slice splitting balances GEMM against conversion, making the
whole $O(p^{1.5})$]
\label{find:phase}
The essential advantage of \OzII{} over \OzI{} is that the number of products
required falls from $O(p^2)$ to $O(p)$ in the target precision.  In the
unsplit configuration $S=1$, however, the price is that the residue
computation and the CRT grow as $O(p^2)$: per element an $O(p/64)$-word
multiple-precision integer is divided by $O(p)$ moduli.  The upper block of
Table~\ref{tab:phase} confirms this.  Doubling $p$ multiplies the GEMM time
by $1.80$--$2.04$ ($O(p)$) but the conversion time by $2.87$--$3.61$
($O(p^2)$), and for $n\le2048$ the GEMM accounts for only $1.6$--$7.8$\% of
the total.  That is, the $n^3$ term drops from $O(p^2n^3)$ to $O(pn^3)$ at
the cost of an $O(p^2n^2)$ conversion term proportional to $n^2$, the two
balancing only at $n\approx6000$ ($p=1024$) to $n\approx31000$ ($p=8192$).

Slice splitting, \eqref{eq:TL}, removes this imbalance.  As the lower block
of Table~\ref{tab:phase} shows, at the optimal slice count $S^{*}$ the two
grow together at $2.2$--$3.1$, that is towards $p^{1.5}$, and the GEMM share
settles at $26$--$35$\%.  This is the direct expression of $T^{*}=2\sqrt{AB}$
being the point at which the two are balanced.  The $n^3$ and $n^2$ terms
become $O(p^{1.5}n^3)$ and $O(p^{1.5}n^2)$ respectively; measured at $n=512$,
$p=8192$ falls from $13.32$~s to $3.89$~s ($3.4\times$), and at $p=16384$ --
a regime the unsplit INT8 path cannot reach at all for lack of capacity --
from $31.28$~s for unsplit binary64 to $10.19$~s ($3.1\times$).
\end{finding}

\begin{table}[htbp]
\centering
\caption{Within the same \OzII{}, a different way of obtaining the required
capacity changes how the number of low-precision GEMMs grows.  Left: the GPU
implementation (\texttt{mpc\_cuda}), which fixes the moduli to the 54 values
below $255$ and splits the significand into slices, issuing the $S(S+1)/2$
digit-group terms of \eqref{eq:gemmcount}.  Right: the CPU implementation
(\texttt{oz2}), which extends the moduli up to $16191$ by two-digit encoding
and performs no slice splitting.  ``Ratio'' is the ratio of counts when $p$ is
doubled.}
\label{tab:gemmcount}
\small
\setlength{\tabcolsep}{5pt}
\begin{tabular}{@{}rrrrrr@{}}
\toprule
\multicolumn{3}{c}{GPU: slice splitting (54 fixed moduli)}
 & \multicolumn{3}{c}{CPU: extended modulus table (no splitting)} \\
\cmidrule(lr){1-3}\cmidrule(l){4-6}
$p$ & \#GEMM & ratio & $p$ & \#GEMM & ratio \\
\midrule
128 & 47 & --- & 1024 & 432 & --- \\
256 & 126 & 2.68 & 2048 & 876 & 2.03 \\
512 & 400 & 3.17 & 4096 & 1782 & 2.03 \\
1024 & 1288 & 3.22 & 8192 & 3687 & 2.07 \\
\bottomrule
\end{tabular}
\end{table}

\begin{finding}[The two routes to capacity are an accuracy-cost trade]
\label{find:gemmcount}
There are two ways to obtain the capacity a target precision requires:
enlarge the modulus table, or split the significand into slices and reduce
the width of one pass.  The growth of the low-precision GEMM count differs
qualitatively.  As Table~\ref{tab:gemmcount} shows, enlarging the modulus
table gives $2.03$--$2.07\times$ per doubling of $p$, that is $O(p)$, whereas
slice splitting gives $2.68$--$3.22\times$, asymptotically $O(p^2)$.  This
follows directly from \eqref{eq:gemmcount} containing $S(S+1)/2$ with
$S\propto p$.

On the other hand the slice-splitting side divides a chunk of width $w$ by a
small number of moduli per element, so its residue computation and CRT are
light, whereas the modulus-table side divides a $p$-bit integer by $O(p)$
moduli and its conversion grows as $O(p^2)$.  The two are therefore a trade
between more GEMMs with lighter conversion and heavier conversion with fewer
GEMMs, and \eqref{eq:TL} gives the balance point.  In this work slice
splitting is implemented on both CPU and GPU and $S$ is chosen by
\eqref{eq:Lhint}, which automates the trade.  The CPU implementation
originally used only modulus-table extension ($S=1$); by
Finding~\ref{find:phase} that was an endpoint far from the balance point.
\end{finding}

\subsubsection{binary64 back-end: taking larger moduli}
\label{sec:fp64}

Opposite to the direction of using a low-precision engine, one may take larger
moduli and reduce their number.  This is precisely the CPU configuration of
the original paper~\cite{ozaki2-2025}: the centered residues are placed
exactly in binary64 and multiplied with an ordinary DGEMM.  The upper bound on
the modulus follows from exactness of the FP64 accumulation,
\begin{equation}
K\left(\frac{m}{2}\right)^2 < 2^{53}
\quad\Longleftrightarrow\quad
m < \sqrt{2^{55}/K} ,
\label{eq:fp64mod}
\end{equation}
so that at $K=512$ we get $m<8.4\times10^6$, about 23 bits per modulus---close
to three times INT8's 8 bits.  The GEMM count consequently falls to a third or
less:

\begin{center}
\begin{tabular}{@{}lrr@{}}
\toprule
Precision & INT8 & FP64 \\
\midrule
DD (106) & 39 & \textbf{13} \\
TD (159) & 81 & \textbf{18} \\
QD (212) & 105 & \textbf{36} \\
\bottomrule
\end{tabular}
\end{center}

At the same time the work of the split and the CRT reconstruction falls in
proportion to the number of moduli, so this can be advantageous where the
overhead of \OzII{} dominates.  For (S4) an off-the-shelf DGEMM is used
directly.  On the CPU we call \texttt{cblas\_dgemm} of
OpenBLAS~\cite{openblas}; since the outer row-block loop is parallelised with
OpenMP, BLAS itself is pinned to serial execution with
\texttt{openblas\_set\_num\_threads(1)} to avoid nested parallelism (the
previous setting is restored on return).  On the GPU we call cuBLASLt with
\texttt{CUDA\_R\_64F}/\texttt{CUBLAS\_COMPUTE\_64F}
(\texttt{src/oz2\_gdtq.cuh}, \texttt{OZG\_FP64}).  An implementation advantage
of this format alone is that it needs no \OzII{}-specific kernel: (S4) is
completed by a stock DGEMM.

\subsection{Performance of LU decomposition}
\label{sec:lubench}

To correspond with \S\ref{sec:gemmbench}, we evaluate the LU decomposition at
the same precision classes (DD/TD/QD) and the same matrix sizes.  To isolate
performance from the influence of the condition number, we use diagonally
dominant random (well-conditioned) matrices here; accuracy and performance on
ill-conditioned matrices are treated in \S\ref{sec:lotkin}.

\subsubsection{CPU: DD/TD/QD}

On the CPU we compare three methods.
\begin{description}
\item[dtq] Native dtq arithmetic for both the panel and the Schur complement
  update; this corresponds to the existing implementation.
\item[\OzII{}-MPFR] The whole matrix is held in MPFR, the panel is done in
  MPFR, and the Schur complement update uses the \OzII{} GEMM.
\item[\OzII{}-dtq] The matrix stays in dtq, the panel uses native dtq (i.e.\
  exactly the same code as the dtq method), and only the Schur complement
  update is delegated to \OzII{}.  Because it uses the direct conversion path
  of \S\ref{sec:directconv}, arrays of \texttt{dd\_real} and the like are
  simply passed as \texttt{double[nc]} arrays and MPFR never appears.
\end{description}

\begin{table}[htbp]
\centering
\caption{CPU: breakdown of the DD/TD/QD LU decomposition at $n=512$
(well-conditioned random matrix).  The set of methods matches
Table~\ref{tab:cpu_lu_qd1024} ($n=1024$).  The three blocked methods use panel
width $b=256$.  Values are from the most recent run on the same matrix on both
machines.  Bold marks the fastest total time for each machine and precision.}
\label{tab:cpu_lu_ddtdqd}
\small
\setlength{\tabcolsep}{5pt}
\begin{tabular}{@{}lllrrrr@{}}
\toprule
Machine & Prec. & Method & Total [s] & panel [s] & Schur [s] & bits \\
\midrule
Arm/GB10 & DD & dtq ($b{=}1$) & 0.036 & 0.009 & 0.027 & 101 \\
 &  & dtq (blk) & \textbf{0.034} & 0.030 & 0.005 & 101 \\
 &  & \OzII{}-MPFR & 0.174 & 0.145 & 0.029 & 101 \\
 &  & \OzII{}-dtq & 0.037 & 0.032 & 0.005 & 102 \\
\cmidrule(l){2-7}
 & TD & dtq ($b{=}1$) & \textbf{0.040} & 0.005 & 0.035 & 156 \\
 &  & dtq (blk) & 0.058 & 0.045 & 0.013 & 156 \\
 &  & \OzII{}-MPFR & 0.188 & 0.158 & 0.030 & 154 \\
 &  & \OzII{}-dtq & 0.057 & 0.049 & 0.007 & 156 \\
\cmidrule(l){2-7}
 & QD & dtq ($b{=}1$) & \textbf{0.068} & 0.005 & 0.062 & 210 \\
 &  & dtq (blk) & 0.115 & 0.091 & 0.023 & 210 \\
 &  & \OzII{}-MPFR & 0.211 & 0.189 & 0.022 & 207 \\
 &  & \OzII{}-dtq & 0.104 & 0.086 & 0.019 & 210 \\
\midrule
x86/H100 & DD & dtq ($b{=}1$) & 0.043 & 0.013 & 0.030 & 101 \\
 &  & dtq (blk) & 0.023 & 0.017 & 0.006 & 101 \\
 &  & \OzII{}-MPFR & 0.126 & 0.114 & 0.012 & 101 \\
 &  & \OzII{}-dtq & \textbf{0.018} & 0.014 & 0.004 & 102 \\
\cmidrule(l){2-7}
 & TD & dtq ($b{=}1$) & 0.039 & 0.007 & 0.032 & 156 \\
 &  & dtq (blk) & 0.028 & 0.018 & 0.010 & 156 \\
 &  & \OzII{}-MPFR & 0.117 & 0.101 & 0.015 & 154 \\
 &  & \OzII{}-dtq & \textbf{0.024} & 0.018 & 0.006 & 156 \\
\cmidrule(l){2-7}
 & QD & dtq ($b{=}1$) & 0.073 & 0.008 & 0.065 & 210 \\
 &  & dtq (blk) & 0.046 & 0.024 & 0.021 & 210 \\
 &  & \OzII{}-MPFR & 0.120 & 0.105 & 0.015 & 207 \\
 &  & \OzII{}-dtq & \textbf{0.032} & 0.024 & 0.007 & 210 \\
\bottomrule
\end{tabular}
\end{table}

\begin{table}[htbp]
\centering
\caption{CPU LU decomposition (well-conditioned random matrices, $n=1024$,
$b=512$), split into the panel factorization and the Schur complement update
that \OzII{} replaces.  dtq (blk) is the blocked dtq native implementation,
\OzII{}-MPFR holds the matrix as \texttt{mpfr\_t}, and \OzII{}-dtq keeps the
expansion format and uses the direct conversion of \S\ref{sec:directconv}.
Bold marks the fastest value for each precision; bits is the number of
correct bits.}
\label{tab:cpu_lu_qd1024}
\small
\setlength{\tabcolsep}{4pt}
\begin{tabular}{@{}lllrrrr@{}}
\toprule
Machine & Prec. & Method & Total [s] & panel [s] & Schur [s] & bits \\
\midrule
Arm/GB10 & DD & dtq ($b{=}1$) & 0.178 & 0.016 & 0.162 & 101 \\
 &  & dtq (blk) & 0.106 & 0.069 & 0.037 & 101 \\
 &  & \OzII{}-MPFR & 1.449 & 1.389 & 0.059 & 101 \\
 &  & \OzII{}-dtq & \textbf{0.102} & 0.079 & 0.023 & 100 \\
\cmidrule(l){2-7}
 & TD & dtq ($b{=}1$) & 0.263 & 0.016 & 0.247 & 155 \\
 &  & dtq (blk) & 0.327 & 0.222 & 0.105 & 155 \\
 &  & \OzII{}-MPFR & 1.654 & 1.585 & 0.069 & 154 \\
 &  & \OzII{}-dtq & \textbf{0.249} & 0.211 & 0.037 & 155 \\
\cmidrule(l){2-7}
 & QD & dtq ($b{=}1$) & 0.502 & 0.022 & 0.479 & 209 \\
 &  & dtq (blk) & 0.633 & 0.458 & 0.176 & 209 \\
 &  & \OzII{}-MPFR & 1.855 & 1.770 & 0.085 & 206 \\
 &  & \OzII{}-dtq & \textbf{0.482} & 0.440 & 0.043 & 209 \\
\midrule
x86/H100 & DD & dtq ($b{=}1$) & 0.115 & 0.021 & 0.095 & 101 \\
 &  & dtq (blk) & 0.064 & 0.038 & 0.026 & 101 \\
 &  & \OzII{}-MPFR & 0.687 & 0.645 & 0.042 & 101 \\
 &  & \OzII{}-dtq & \textbf{0.051} & 0.038 & 0.012 & 100 \\
\cmidrule(l){2-7}
 & TD & dtq ($b{=}1$) & 0.240 & 0.021 & 0.218 & 155 \\
 &  & dtq (blk) & 0.164 & 0.075 & 0.089 & 155 \\
 &  & \OzII{}-MPFR & 0.758 & 0.706 & 0.052 & 154 \\
 &  & \OzII{}-dtq & \textbf{0.094} & 0.075 & 0.019 & 155 \\
\cmidrule(l){2-7}
 & QD & dtq ($b{=}1$) & 0.500 & 0.023 & 0.476 & 209 \\
 &  & dtq (blk) & 0.348 & 0.127 & 0.221 & 209 \\
 &  & \OzII{}-MPFR & 0.791 & 0.736 & 0.055 & 206 \\
 &  & \OzII{}-dtq & \textbf{0.148} & 0.123 & 0.024 & 209 \\
\bottomrule
\end{tabular}
\end{table}

\begin{finding}[The machine matters on the CPU too; \OzII{} clearly wins with
AMX-INT8]
\label{find:cpu_lu_ddtdqd}
As the last column of Table~\ref{tab:cpu_lu_qd1024} ($n=1024$) shows, the
speedup of \OzII{} over the best existing implementation differs greatly
between the machines:
\begin{center}
Arm/GB10 (SVE2 i8mm): $0.78/1.11/1.06\times$ \qquad
x86/H100 (AMX-INT8): $1.32/1.81/2.41\times$ .
\end{center}
On the x86 side, AMX-INT8 is more powerful than SVE2 i8mm, so \OzII{} wins
clearly at every precision.  On Arm/GB10, by contrast, it loses at DD
($0.78\times$) and only barely wins at TD/QD.  Common to both machines is that
the ratio increases with precision, consistent with the trend of the GEMM
alone in \S\ref{sec:gemmbench} (Finding~\ref{find:gemm_cpu}).

The margin is small on Arm/GB10 because even on the \OzII{} side the panel
factorisation (DD 0.090~s,
TD 0.212~s, QD 0.440~s) accounts for 80--90\% of the total, so no matter how
fast the Schur side becomes, this is the ceiling.  The unblocked LU in
particular has almost no panel cost (0.023--0.026~s), and that is what tells.
\end{finding}

\begin{finding}[In multiple precision the unblocked LU is a strong comparison
target]
\label{find:unblocked}
As Table~\ref{tab:cpu_lu_qd1024} shows, at TD/QD the unblocked ($b=1$) dtq LU
is faster than the blocked dtq LU (TD: 0.277~s against 0.361~s; QD: 0.518~s
against 0.594~s).  The main purpose of blocking is cache reuse, but in
multiple precision the cost of one operation is high and the workload has
shifted from memory-bound towards compute-bound, so that benefit does not
outweigh the panel-factorisation cost $O(n^2b)$.  The tendency is stronger for
smaller $n$: at $n=512$ on Arm/GB10 at QD in
Table~\ref{tab:cpu_lu_ddtdqd}, the unblocked dtq takes 0.070~s, more than
twice as fast as \OzII{}-dtq's 0.144~s, because the unblocked version spends
only $0.006$~s on the panel whereas \OzII{}-dtq needs $0.135$~s.  Any
performance evaluation of multiple-precision LU must therefore include the
unblocked version as a comparison target.  Had we compared only against the
blocked naive version, as we did initially, we would have overestimated the
advantage of \OzII{}.
\end{finding}

\begin{finding}[The design of the conversion path decides viability]
\label{find:cpu_lu_schur}
Within the same \OzII{}, the MPFR-mediated \OzII{}-MPFR is 2.7--15 times
slower than dtq in total time (Table~\ref{tab:cpu_lu_qd1024}).  This is
because the panel runs in MPFR and because every Schur complement update
performs per-element \texttt{mpfr\_init2}/\texttt{mpfr\_clear} on the three
matrices $L_{21},U_{12},U$ (at $n=1024$, $b=512$ the output alone is
$262{,}144$ elements).  The direct conversion path eliminates this heap
allocation entirely, shortening the same Schur complement update from 0.056~s
to 0.024~s at DD and from 0.078~s to 0.042~s at QD.  That the success of
\OzII{} at multi-component precision is decided not by the algorithm but by
how the data representations are connected is the most practical lesson of
this work.
\end{finding}

\subsubsection{GPU: gdtq (multi-component)}

The GPU side corresponds exactly to the CPU side, using gdtq's multi-component
types directly.  The two methods compared are
\begin{description}
\item[gdtq] Native gdtq arithmetic for both the panel and the Schur complement
  update (using \texttt{operator*} and the like in device kernels); this
  corresponds to the existing implementation.
\item[\OzII{}] The same panel code, with only the Schur complement update
  replaced by the \OzII{} of \S\ref{sec:impl_gpu} (direct conversion from the
  expansion format).
\end{description}
These correspond one-to-one with the dtq/\OzII{}-dtq pair on the CPU.  The
precisions and matrix sizes are identical to the CPU case as well
(DD/TD/QD $=106/159/212$ bits, $n=512,1024$).

\begin{table}[htbp]
\centering
\caption{GPU: LU decomposition with gdtq at DD/TD/QD, $n=1024$, $b=512$
(well-conditioned random matrix).  The panel code is the same for all
methods.  ``Schur ratio'' is the speedup of the Schur complement update over
gdtq.  Bold marks the fastest Schur complement update for each machine and
precision.  The number of correct bits agrees across methods and is omitted
(DD 100.6 / QD 208.7 bits; only for TD does gdtq differ slightly, 154.9
against 155.0 for the \OzII{} variants).  All methods were measured after a
warm-up.  Note that the ranking reverses between the two machines
(Finding~\ref{find:gpu_reversal}).}
\label{tab:gpu_lu_ddtdqd}
\small
\setlength{\tabcolsep}{5pt}
\begin{tabular}{@{}lllrrrr@{}}
\toprule
Machine & Prec. & Method & Total [s] & panel [s] & Schur [ms] & ratio \\
\midrule
Arm/GB10 & DD & gdtq          & 0.1242 & 0.1007 & 23.47 & --- \\
         &    & \OzII{}/INT8  & 0.1036 & 0.1004 & \textbf{3.18} & 7.4 \\
         &    & \OzII{}/FP16  & 0.1039 & 0.1005 & 3.32 & 7.1 \\
         &    & \OzII{}/FP8   & 0.1046 & 0.1003 & 4.24 & 5.5 \\
         &    & \OzII{}/FP64  & 0.1117 & 0.1002 & 11.50 & 2.0 \\
\cmidrule(l){2-7}
         & TD & gdtq          & 0.3012 & 0.2162 & 85.00 & --- \\
         &    & \OzII{}/INT8  & 0.2204 & 0.2158 & \textbf{4.63} & 18.4 \\
         &    & \OzII{}/FP16  & 0.2207 & 0.2160 & 4.71 & 18.1 \\
         &    & \OzII{}/FP8   & 0.2229 & 0.2161 & 6.70 & 12.7 \\
         &    & \OzII{}/FP64  & 0.2399 & 0.2170 & 22.90 & 3.7 \\
\cmidrule(l){2-7}
         & QD & gdtq          & 0.5807 & 0.3949 & 185.80 & --- \\
         &    & \OzII{}/INT8  & 0.4014 & 0.3951 & \textbf{6.34} & 29.3 \\
         &    & \OzII{}/FP16  & 0.4014 & 0.3949 & 6.51 & 28.5 \\
         &    & \OzII{}/FP8   & 0.4041 & 0.3950 & 9.13 & 20.4 \\
         &    & \OzII{}/FP64  & 0.4256 & 0.3951 & 30.50 & 6.1 \\
\midrule
x86/H100 & DD & gdtq          & 0.0480 & 0.0474 & \textbf{0.48} & --- \\
         &    & \OzII{}/INT8  & 0.0489 & 0.0474 & 1.48 & 0.33 \\
         &    & \OzII{}/FP16  & 0.0492 & 0.0475 & 1.68 & 0.29 \\
         &    & \OzII{}/FP8   & 0.0502 & 0.0476 & 2.56 & 0.19 \\
         &    & \OzII{}/FP64  & 0.0485 & 0.0474 & 1.06 & 0.46 \\
\cmidrule(l){2-7}
         & TD & gdtq          & 0.0512 & 0.0497 & \textbf{1.42} & --- \\
         &    & \OzII{}/INT8  & 0.0522 & 0.0501 & 2.12 & 0.67 \\
         &    & \OzII{}/FP16  & 0.0521 & 0.0499 & 2.25 & 0.63 \\
         &    & \OzII{}/FP8   & 0.0538 & 0.0499 & 3.95 & 0.36 \\
         &    & \OzII{}/FP64  & 0.0513 & 0.0498 & 1.46 & 0.97 \\
\cmidrule(l){2-7}
         & QD & gdtq          & 0.0562 & 0.0533 & 2.84 & --- \\
         &    & \OzII{}/INT8  & 0.0565 & 0.0537 & 2.78 & 1.02 \\
         &    & \OzII{}/FP16  & 0.0563 & 0.0533 & 2.95 & 0.96 \\
         &    & \OzII{}/FP8   & 0.0586 & 0.0535 & 5.08 & 0.56 \\
         &    & \OzII{}/FP64  & 0.0553 & 0.0534 & \textbf{1.83} & 1.55 \\
\bottomrule
\end{tabular}
\end{table}

\begin{finding}[On the GPU the ranking reverses between machines]
\label{find:gpu_reversal}
The most important fact shown by Table~\ref{tab:gpu_lu_ddtdqd} is that the
conclusion is exactly opposite on the two machines.
\begin{itemize}
\item Arm/GB10: the \OzII{}/INT8 Schur complement update is $7.4/17.9/29.4$
  times faster than gdtq native (DD/TD/QD), and $1.20/1.36/1.45$ times faster
  in total.
\item x86/H100: gdtq native is faster, at $0.34\times$ for DD and
  $0.71\times$ for TD; only at QD do they become comparable at $1.06\times$.
\end{itemize}
The cause is the difference in the strength of the FP64/FP32 units.  gdtq's
DD/TD/QD arithmetic runs on the FP64 units, so it is slow on Arm/GB10 with its
weak FP64 and very fast on x86/H100 with its strong FP64.  \OzII{}, by
contrast, uses tensor cores on both machines, leaving relatively less headroom.
In other words, whether \OzII{} wins is decided by the ratio of low-precision
to high-precision engine performance on the target machine.

In absolute time, x86/H100 is 2--10 times faster than Arm/GB10, finishing the
whole $n=1024$ LU in 0.05~s.  At that scale the Schur complement update is a
small problem of about $512^3$, and the overhead of \OzII{}'s pre-processing
and CRT reconstruction weighs relatively more.  We expect both machines to
swing towards \OzII{} at larger $n$, but this was not verified within the
range measured here.
\end{finding}

\begin{finding}[The fastest back-end changes completely with the machine]
\label{find:fp64}
Looking at the Schur complement update of Table~\ref{tab:gpu_lu_ddtdqd} per
machine, the fastest back-end is exactly opposite.
\begin{center}
Arm/GB10: INT8 fastest (FP64 is 3.6--4.8$\times$ slower)\qquad
x86/H100: FP64 fastest (at QD, $1.52\times$ INT8 and $1.55\times$ gdtq)
\end{center}
FP64 underperforms on Arm/GB10 because that machine's FP64 performance is of
the order of $0.5$~TFLOPS, cancelling the benefit of reducing the GEMM count
to one third (\S\ref{sec:fp64}).  On x86/H100, FP64 reaches tens of TFLOPS
through the tensor cores, two orders of magnitude above Arm/GB10.  On that
machine the overhead of \OzII{} (split and CRT reconstruction) dominates, as
in Finding~\ref{find:gpu_reversal}, and the binary64 back-end cuts precisely
that dominant term by reducing the number of moduli to one third.

As a result, at QD on x86/H100, \OzII{}/FP64 is fastest at 1.83~ms, beating
both gdtq native (2.84~ms) and INT8 (2.78~ms), and it also overtakes gdtq in
total time, 0.0553~s against 0.0562~s.  At TD, FP64's 1.46~ms is $1.45\times$
better than INT8's 2.12~ms and essentially matches gdtq (1.42~ms).  Only at DD
is gdtq native fastest.

Using a low-precision engine is therefore not the essence of \OzII{}.  The
essence is to choose, according to the target precision and the balance of
engines on the machine, the format that maximises the product of
bits-per-modulus and throughput; on machines with strong FP64 the right answer
is to take larger moduli and reduce their number.
\end{finding}

\begin{finding}[FP8 cannot win in \OzII{} even where its raw performance is
highest]
\label{find:fp8}
The newly implemented FP16 and FP8 back-ends (\S\ref{sec:fp816}) give exactly
the same number of correct bits as INT8 at all of DD/TD/QD, confirming the
correctness of the implementation.  Their speeds, from
Table~\ref{tab:gpu_lu_ddtdqd}, are
\begin{center}
INT8 $\simeq$ FP16 $<$ FP8 (Schur complement update; FP8 is
$1.35$--$1.44\times$ INT8),
\end{center}
and all three beat gdtq native by $5.5$--$29.5\times$.  FP8 is slightly worse
because the two-digit decomposition of \eqref{eq:fp8digits} needs four
sub-GEMMs per modulus; given that FP8's effective performance on Arm/GB10 is
only $1.29\times$ INT8, most of that fourfold gap is in fact recovered by
engine performance.

This point should be tested precisely on the machine where FP8 is strongest.
The raw GEMM performance of x86/H100 is 1525.7~TFLOP/s for FP8, $1.94\times$
FP16 (786.2~TFLOP/s) and the fastest of all formats
(Table~\ref{tab:probe}).  Yet even there, the \OzII{} Schur complement update
is, from Table~\ref{tab:gpu_lu_ddtdqd},
\begin{center}
DD: INT8 1.41~ms / FP16 1.56~ms / FP8 2.14~ms\qquad
QD: INT8 2.68~ms / FP16 2.64~ms / FP8 4.16~ms,
\end{center}
so FP8 is the slowest.  The fourfold sub-GEMM count from the two-digit
decomposition of \eqref{eq:fp8digits} cannot be recovered by a $1.94\times$
throughput advantage.

We therefore conclude that FP8's disadvantage in \OzII{} is not a matter of
engine performance but is structural, arising from a representational capacity
of only three significand bits.  This does not contradict
\cite{ozfp8-2026}: the binary64 target of that work ($p=53$) needs the fewest
digit groups, and its comparison target is a native FP64 GEMM rather than an
INT8 version of \OzII{}.  At DD/TD/QD and above, and where INT8 is available,
INT8 should be chosen.
\end{finding}

\subsection{LU decomposition of the Lotkin matrix and the relative error of
the numerical solution}
\label{sec:lotkin}

Finally we evaluate arbitrary-precision LU decomposition on an
ill-conditioned matrix of practical importance.  The Lotkin matrix has all
ones in its first row and coincides with the Hilbert matrix from the second
row on ($A_{ij}=1/(i+j+1)$ for $i\ge1$, 0-indexed), and its condition number
grows rapidly, $\log_2\mathrm{cond}(A)\approx5.1n$.  At the fixed precisions
of DD/TD/QD no solution at all is obtained for $n\gtrsim40$, so arbitrary
precision is essential.  With the true solution $x=(1,\ldots,1)^\top$, we
generate $b=Ax$ at high working precision, round it to the target precision,
solve for $\hat x$ by $PA=LU$ with forward and backward substitution, and take
$\min_i(-\log_2\abs{\hat x_i-1})$ as the ``number of correct bits''.

\begin{table}[htbp]
\centering
\caption{LU decomposition of the Lotkin matrix (total time [s]).  The
precision is $p \approx 1.2\log_2\mathrm{cond}(A) \approx 1.2\times5.1n$
rounded to a multiple of 64.  ``naive $b{=}1$'' is unblocked; ``naive blk''
shares the panel with the \OzII{} version and replaces only the Schur
complement update by a naive inner product (\S\ref{sec:impl_gpu}).  ``BE'' is
the back-end \OzII{} used, chosen per machine by measurement
(\S\ref{sec:int8limit}).  ``best/\OzII{}'' is the speed-up over the faster of
the two naive variants.  Bold marks the fastest value in each row; bits is
the number of correct bits and the relative error is
$\varepsilon = 2^{-\text{bits}}$.  The breakdown into panel factorization and
Schur complement update is shown in Figure~\ref{fig:lotkin}.}
\label{tab:lotkin}
\small
\setlength{\tabcolsep}{4pt}
\begin{tabular}{@{}rrlrrrlrr@{}}
\toprule
$n$ & $p$ & Machine & \multicolumn{2}{c}{naive} & \OzII{} & BE & best & bits \\
\cmidrule(lr){4-5}
 & & & $b{=}1$ & blk & & & /\OzII{} & \\
\midrule
512 & 3136 & CPU GB10 & 3.032 & 3.507 & \textbf{2.258} & INT8 & \textbf{1.34} & 546 \\
 & & CPU H100 & 1.853 & 1.632 & \textbf{0.8529} & b64 & \textbf{1.91} & 546 \\
 & & GPU GB10 & 2.171 & 2.08 & \textbf{1.651} & INT8 & \textbf{1.26} & 545 \\
 & & GPU H100 & 1.323 & 1.345 & \textbf{1.308} & INT8 & \textbf{1.01} & 545 \\
\midrule
1024 & 6272 & CPU GB10 & 59.96 & 65.67 & \textbf{33.39} & INT8 & \textbf{1.80} & 1083 \\
 & & CPU H100 & 31.82 & 30.17 & \textbf{14.91} & b64 & \textbf{2.02} & 1083 \\
 & & GPU GB10 & 381.2 & 376.8 & \textbf{142.3} & INT8 & \textbf{2.65} & 1077 \\
 & & GPU H100 & 30.54 & 32.32 & \textbf{18.68} & INT8 & \textbf{1.63} & 1077 \\
\midrule
2048 & 12544 & CPU GB10 & 1348 & 1432 & \textbf{619} & INT8 & \textbf{2.18} & 2142 \\
 & & CPU H100 & 734.5 & 661.8 & \textbf{250.6} & INT8 & \textbf{2.64} & 2142 \\
 & & GPU GB10 & 17604 & 17649 & \textbf{6206} & INT8 & \textbf{2.84} & 2142 \\
 & & GPU H100 & 919.3 & 944.1 & \textbf{445.3} & INT8 & \textbf{2.06} & 2142 \\
\bottomrule
\end{tabular}
\end{table}

\begin{figure}[htbp]
\centering
\includegraphics[width=\linewidth]{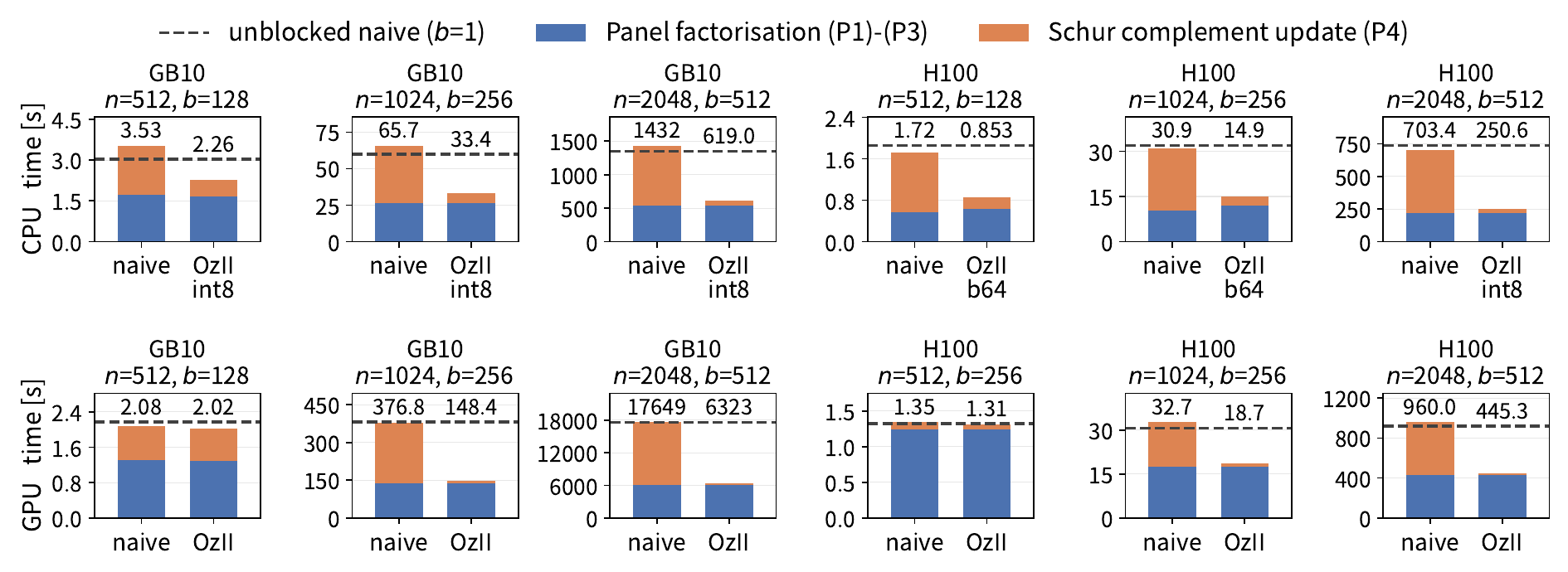}
\caption{Breakdown of the LU decomposition of the Lotkin matrix.  The upper
row is the CPU and the lower the GPU.  Each pair is the blocked naive and
\OzII{} at the same panel width $b$; the panel factorization (P1)--(P3) is
identical code in both.  The dashed line is the total time of the unblocked
naive ($b=1$), which is not drawn as a bar because it does not report the
breakdown separately.  The numbers above the bars are total times in seconds,
with the chosen back-end noted for \OzII{}.  The vertical axis is linear;
because the values differ by three orders of magnitude across $n$, each panel
has its own scale.  $b$ is chosen to minimise the \OzII{} total and the naive
bar uses the same $b$.
$b$ is chosen per environment to minimise the \OzII{} total, so it may
differ between CPU and GPU (the GPU Lotkin runs cover only $b=256$ for
$n=512,1024$ and $b=512$ for $n=2048$).  Where a configuration was run
repeatedly, the run with the smallest Schur complement update is taken as
representative.  \OzII{} replaces only (P4), which shrinks markedly
while the total remains bound by the panel factorization -- a pattern common
to every condition.}
\label{fig:lotkin}
\end{figure}

\begin{finding}[Accuracy agrees between CPU and GPU; larger $n$ favours
\OzII{}]
\label{find:lotkin}
As Table~\ref{tab:lotkin} shows, the number of correct bits agrees between CPU
and GPU to within 1--3 bits, and in every case slightly exceeds the
theoretical bound $p-\log_2\mathrm{cond}$ (525, 1050 and 2099 respectively), a
reasonable value.  The methods agree with one another as well, so the speedup
costs no accuracy.  In terms of relative error, $n=2048$ reaches $10^{-645}$,
a regime that binary64 cannot even represent.

Larger $n$ favours \OzII{}.  On the CPU the speedup grows
$1.34\to1.80\to2.18\times$ on Arm/GB10 and $1.91\to2.02\to2.64\times$ on
x86/H100; on the GPU it grows $1.01\to1.63\to2.06\times$ on x86/H100 and
$1.26\to2.65\to2.84\times$ on Arm/GB10.  All four environments are
monotone in $n$; the largest value in the table, $2.84\times$, is obtained on the
Arm/GB10 GPU at $n=2048$.  This is a direct consequence of the structure in
which the panel factorisation scales as $O(n^2b)$ and the Schur complement
update as $O(n^3)$.
\end{finding}

\begin{remark}[The numbers in this subsection were initially wrong]
An early version reported 4096 correct bits on the CPU (the working-precision
ceiling).  The error $\approx2^{-1485}$ was being converted to binary64 with
\texttt{mpfr\_get\_d()} before taking $\log_2$, so it fell below binary64's
lower limit $2^{-1074}$, collapsed to 0, and was misjudged as an exact
solution.  After correcting this to use the MPFR exponent
(\texttt{mpfr\_get\_exp}) directly, the CPU gave 1506 bits and the GPU 1504,
in agreement.
\end{remark}

\begin{finding}[Larger and higher-precision problems favour it further]
\label{find:lotkin_large}
Looking at the breakdown of the Schur complement update in
Table~\ref{tab:lotkin}, at $n=2048$ on the H100 CPU the blocked naive at the same $b=512$ takes
485.6~s against 32.6~s for \OzII{} ($14.9\times$).  The total stays at
$2.64\times$ because the panel factorisation, at 218.0~s, accounts for 87\% of
the \OzII{} total---the structure of Finding~\ref{find:profile} holding even at the
largest size.  The GPU shows the same structure: on Arm/GB10 at $n=2048$ the
Schur complement update shrinks from the naive 11555.1~s to 234.1~s (GEMM
48.8~s, CRT 183.5~s, conversion 1.8~s), a factor of $49.4$, while the panel
factorisation at 6088.0~s accounts for 98\% of the \OzII{} total, so the
overall speedup remains $2.84\times$.  The $2.06\times$ on x86/H100 becomes
$2.12\times$ if the comparison target is replaced by the blocked naive
(944.1~s), confirming that it is not an apparent speedup caused by blocking
(Finding~\ref{find:gpu_caveat}).
\end{finding}

\begin{finding}[On Arm/GB10 the GPU is slower than the CPU]
\label{find:gb10_gpu_slow}
Comparing $n=2048$ across machines in Table~\ref{tab:lotkin}, x86/H100 is
roughly balanced---919.3~s on the GPU against 734.5~s on the CPU---whereas on
Arm/GB10 the GPU takes 17604.1~s against the CPU's 1348.4~s, $13.1$ times
slower.  Between the GPUs, Arm/GB10 needs $19.2\times$ the time of x86/H100,
while between the CPUs the ratio is only $1.84\times$.  This is because
Arm/GB10's weak binary64 capability (Table~\ref{tab:env}) shows up directly in
the panel factorisation, which cannot be moved to \OzII{}.  Indeed the \OzII{}
Schur complement update, using the INT8 path, holds up well even on Arm/GB10
at 234.6~s; the bottleneck is entirely the 6088.0~s of panel factorisation.
The structure of Finding~\ref{find:panel}---that the panel factorisation is
the common bottleneck---is thus most visible on machines with weak binary64,
suggesting that this configuration is where lowering the precision of the
panel factorisation itself, discussed in \S\ref{sec:future}, would help most.
\end{finding}

\begin{finding}[The GPU speedup is not an effect of blocking]
\label{find:gpu_caveat}
Initially the only GPU comparison target was the unblocked implementation
(\texttt{lu\_normal} of \texttt{mpc\_cuda\_lu}), so the speedup mixed ``the
effect of blocking itself'' with ``the effect of doing the Schur complement
update with \OzII{}''.  We therefore also implemented on the GPU a ``blocked
naive'' that shares the panel and pivoting code completely with the \OzII{}
version and replaces only the Schur complement update by a naive
\texttt{cu\_freal} inner product (\S\ref{sec:impl_gpu}), performing the same
separation as in the CPU experiments.

The result is in the GPU rows of Table~\ref{tab:lotkin}: at all six points the
difference between unblocked and blocked naive lies between $-4.2$\% and
$+5.8$\%.  On Arm/GB10, $n=512$ is 2.17~s against 2.08~s, $n=1024$ is 381.2~s
against 376.8~s, and $n=2048$ is 17604.1~s against 17649.1~s; on x86/H100,
$n=512$ is 1.32~s against 1.35~s, $n=1024$ is 30.5~s against 32.3~s, and
$n=2048$ is 919.3~s against 944.1~s.  The improvement in locality from
blocking is roughly cancelled by the overhead of splitting the Schur
complement update by panel width.  The \OzII{} speedups ($1.26$--$2.84\times$
on Arm/GB10, $1.01$--$2.06\times$ on x86/H100) are therefore almost entirely
due to the method of the Schur complement update, with a negligible
contribution from blocking itself.  The DD/TD/QD comparison of
Table~\ref{tab:gpu_lu_ddtdqd} was obtained with the same separation.
\end{finding}

\begin{finding}[The panel factorisation is the common bottleneck]
\label{find:panel}
Across Tables~\ref{tab:cpu_lu_ddtdqd}--\ref{tab:lotkin}, the panel
factorisation time is nearly identical between naive and \OzII{} (naturally,
since the code is the same), but it takes a large share of the total (77\% for
the CPU Lotkin case, 96\% for the GPU at $n=1024$/$p=8192$).  However much the
Schur complement update is accelerated, the overall speedup is limited wherever
this part is the bottleneck.  In our CPU implementation the panel
factorisation was initially single-threaded and accounted for more than 80\%
of the total; OpenMP parallelisation along rows shortened it by about eight
times (all numbers in this section are after that parallelisation).
\end{finding}

\subsection{Summary: profiling the LU decomposition}
\label{sec:profile}

\begin{figure}[H]
\centering
\includegraphics[width=\textwidth]{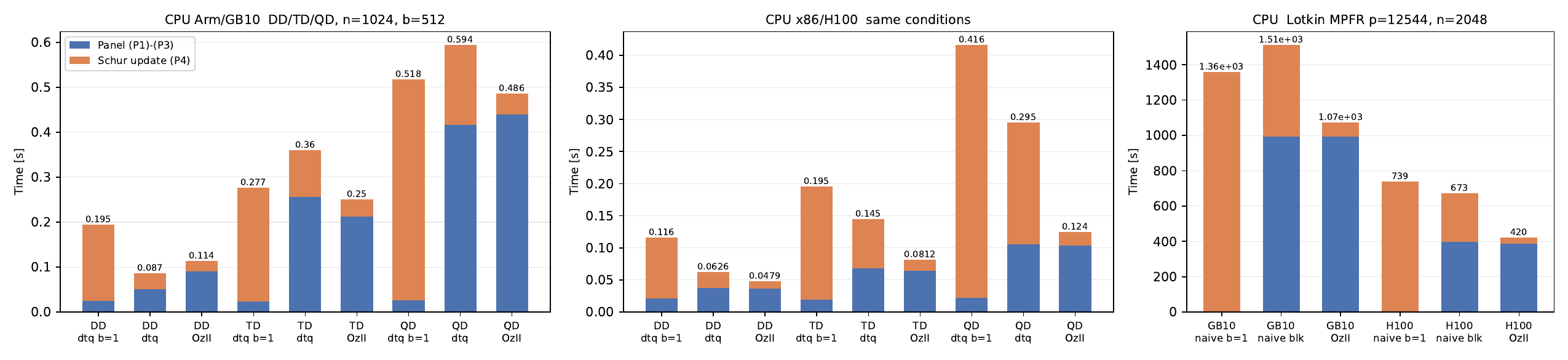}
\caption{Breakdown of the LU decomposition on the CPU, $n=1024$, $b=512$
(DD/TD/QD, well-conditioned random matrix; \texttt{b=1} is unblocked).  Left:
Arm/GB10; centre: x86/H100 under the same conditions; right: the
arbitrary-precision Lotkin case ($p=12544$, $n=2048$) for both machines.  The
lower part of each bar is the panel factorisation (P1)--(P3), the upper part
the Schur complement update (P4); the number is the total time [s].  Note that
the vertical scale differs between panels.  The pattern---the unblocked
(\texttt{b=1}) case has an almost invisible panel but a large Schur part,
and blocking reverses this---is common to all precisions and machines.}
\label{fig:prof_cpu}
\end{figure}

\begin{figure}[H]
\centering
\includegraphics[width=\textwidth]{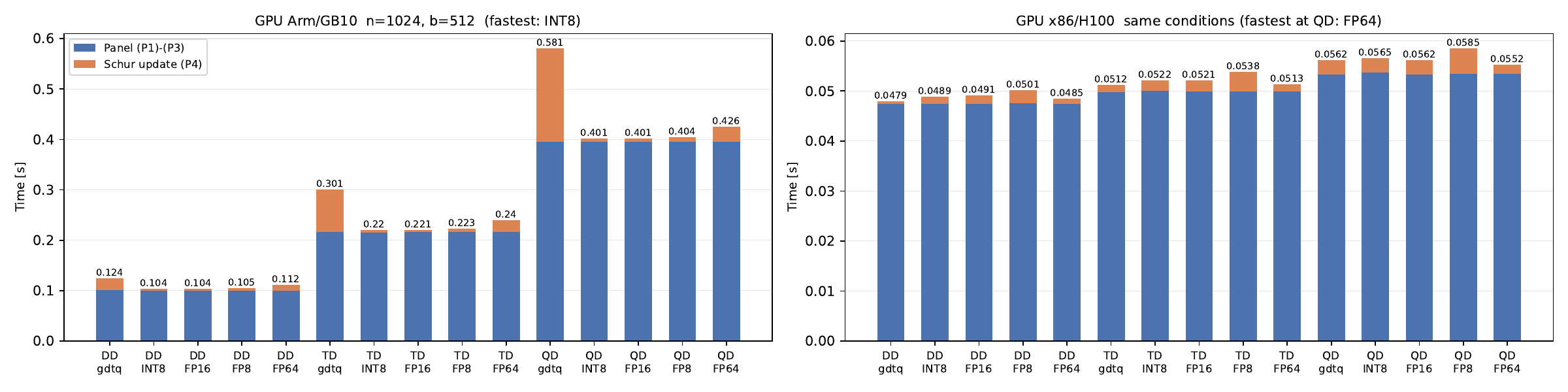}
\caption{Breakdown of the LU decomposition on the GPU, $n=1024$, $b=512$
(DD/TD/QD with gdtq).  Left: Arm/GB10; right: x86/H100 under the same
conditions.  The lower part of each bar is the panel factorisation, the upper
part the Schur complement update.  Besides the vertical scales differing by an
order of magnitude, the gdtq Schur part dominates on Arm/GB10 where INT8 is
fastest, whereas on x86/H100 the gdtq Schur part is barely visible and FP64 is
fastest at QD (Findings~\ref{find:gpu_reversal} and~\ref{find:fp64}).}
\label{fig:prof_gpu}
\end{figure}

To close this section, we present the measurements given individually in
\S\ref{sec:lubench}--\ref{sec:lotkin} as a profile that divides the LU
decomposition into the panel factorisation ((P1)--(P3), identical code in both
methods) and the Schur complement update ((P4), the only part that differs),
shown as stacked bar charts in Figs.~\ref{fig:prof_cpu}
and~\ref{fig:prof_gpu}.  DD/TD/QD use well-conditioned random matrices and
MPFR uses the Lotkin matrix.

\begin{finding}[The Schur complement update all but vanishes, leaving the
panel factorisation]
\label{find:profile}
The structure visible in Figs.~\ref{fig:prof_cpu} and~\ref{fig:prof_gpu} is
clear.
\begin{enumerate}
\item \OzII{} accelerates the Schur complement update by $1.50$--$3.87\times$
  on the CPU and $7.3$--$29.5\times$ on the GPU, and the factor grows with
  precision.
\item As a result, on the \OzII{} side the panel factorisation grows to
  80--98\% of the total: the Schur complement update is no longer the
  bottleneck.
\item The improvement in total time ($0.78$--$1.11\times$ on the CPU,
  $1.20$--$1.45\times$ on the GPU) is far smaller than the factor for the
  Schur part alone precisely because of this Amdahl-type saturation.
\item The bars for the unblocked case ($b=1$) have an almost invisible panel,
  showing visually that the cost of blocking is not negligible in multiple
  precision (Finding~\ref{find:unblocked}).
\end{enumerate}
At QD on the GPU the Schur complement update has shrunk to 1.6\% of the total
($0.0063/0.4020$), and for MPFR to 0.7\% ($3.7/500.8$); further optimisation
of the Schur side would change almost nothing overall.  The profile shows
quantitatively that the panel factorisation is what to tackle next.

This bears directly on the choice of panel width $b$.  For MPFR on the GPU
($n=1024$, $p=8192$), $b=512$ gives a total of 500.8~s (panel 497.1~s), but
lowering it to $b=256$ improves this substantially to 300.4~s (panel 288.4~s,
Schur 11.9~s).  Because the panel factorisation scales as $O(nb^2)$, a smaller
$b$ is faster wherever the panel dominates.  Conversely, where the Schur
complement update dominates, a larger $b$ is advantageous because the
pre-processing and CRT reconstruction cost of \OzII{} does not depend on the
inner dimension $b$.  The optimal $b$ is set by the balance between the two,
and in our measurements this balance point moves with precision, size and
platform.
\end{finding}

\section{Conclusions and future work}
\label{sec:conclusion}

\subsection{Conclusions}

We implemented a multiple-precision LU decomposition based on \OzII{} on both
CPU and GPU, released as the open-source library mpoz2~\cite{mpoz2}, and evaluated CPU and GPU at identical precision and identical
size for both multi-component (DD/TD/QD) and arbitrary precision.  Our
conclusions are as follows.

\begin{enumerate}
\item For matrix multiplication we obtained speedups of
  $3.69/7.01/12.42\times$ for DD/TD/QD at $n=1024$ on the CPU (Arm/GB10) and
  $9.5$--$15.9\times$ (same machine) and $12.8$--$32.2\times$ (x86/H100) for
  MPFR arbitrary precision (Finding~\ref{find:gemm_cpu}), and, with INT8 on the
  GPU, $40$--$66\times$ on Arm/GB10 and $17$--$20\times$ on x86/H100
  (Table~\ref{tab:gpu_gemm}).  The relative error is consistently five to
  seven orders of magnitude smaller than that of the existing implementation.
\item The Schur complement update of the LU decomposition was accelerated at
  $n=1024$ by $1.59/2.80/4.12\times$ on the CPU (dtq, Arm/GB10) and
  $2.12/4.70/9.16\times$ (x86/H100), and by $7.4/17.9/29.4\times$ on the GPU
  (gdtq) for DD/TD/QD on Arm/GB10.  Both use the same
  multi-component types as the existing implementation (dtq on the CPU, gdtq
  on the GPU) directly, without \texttt{cu\_freal} or MPFR
  (\S\ref{sec:directconv}).  The MPFR-mediated implementation on the CPU is,
  conversely, much worse, so the design of the conversion path decides the
  viability of \OzII{} itself (Finding~\ref{find:cpu_lu_schur}).
\item The improvement in \emph{total} LU time, however, depends strongly on
  the machine.  Against the best existing implementation (in many cases the
  unblocked $b=1$), the CPU gives $1.04$--$1.06\times$ on Arm/GB10 and
  $1.27$--$2.36\times$ on x86/H100 at $n=1024$
  (Finding~\ref{find:cpu_lu_ddtdqd}).  On the
  GPU the ranking even reverses: $1.20$--$1.45\times$ on Arm/GB10 becomes
  $0.34$--$1.07\times$ on x86/H100 (Finding~\ref{find:gpu_reversal}).  The
  reason is that gdtq's DD/TD/QD arithmetic runs on the FP64 units while
  \OzII{} uses tensor cores, so the ratio of low- to high-precision engine
  performance decides the outcome.  Moreover, as the profile
  (\S\ref{sec:profile}) shows, after applying \OzII{} the panel factorisation
  accounts for 80--98\% of the total and the Schur complement update is no
  longer the bottleneck (Finding~\ref{find:profile}).  Any evaluation of
  multiple-precision LU must include the unblocked version as a comparison
  target (Finding~\ref{find:unblocked}).
\item The fastest back-end changes with the machine.  In arbitrary
  precision, Arm/GB10 favours INT8 at every $p$ (ratio $0.72$--$0.88$) while
  x86/H100 favours binary64 at every $p$ (ratio $1.24$--$2.19$); on the GPU
  the difference between machines is far more extreme ($0.07$--$0.09$ against
  $1.73$--$1.88$).  The dividing line is the throughput ratio between the
  low-precision engine and binary64 ($13.1\times$ on GB10, $7.5\times$ on
  H100).  The essence of \OzII{} is not ``use a low-precision engine'' but
  ``choose the format that maximises the product of bits-per-modulus and
  throughput''.  We attempted to automate this choice at run time but
  abandoned it: branching on the architecture, the analytic cost model and a
  start-up probe all fail to generalise to other environments
  (\S\ref{sec:int8limit}).
\item We newly implemented FP16 and FP8 back-ends and confirmed that both give
  the same number of correct bits as INT8 (Finding~\ref{find:fp8}).  FP8
  secures the same 362.8-bit CRT capacity as INT8 through the balanced
  base-17 two-digit decomposition of \eqref{eq:fp8digits}.  The speed order is
  INT8 $\simeq$ FP16 $<$ FP8 on both machines.  Notably, even on the machine
  where FP8 has the highest raw GEMM performance of any format (1525.7~TFLOP/s
  on x86/H100, $1.94\times$ FP16), FP8 is the slowest within \OzII{}: a
  $1.94\times$ throughput advantage cannot recover the fourfold sub-GEMM count
  of the two-digit decomposition.  FP8's disadvantage in \OzII{} is structural,
  arising from the representational capacity of three significand bits rather
  than from engine performance.
\item We evaluated the arbitrary-precision LU decomposition of the Lotkin
  matrix at $n=512/1024/2048$ with
  $p\approx1.2\log_2\mathrm{cond}$.  The number of correct bits agrees between
  CPU and GPU to within 1--3 bits and is consistent with the theoretical
  bound.  The relative error reaches $10^{-645}$ at $n=2048$.  Larger $n$
  favours \OzII{}: $1.34\to1.80\to2.18\times$ on the Arm/GB10 CPU,
  $1.91\to2.02\to2.64\times$ on the x86/H100 CPU,
  $1.26\to2.65\to2.84\times$ on the Arm/GB10 GPU, and
  $1.01\to1.63\to2.06\times$ on the x86/H100 GPU
  (Finding~\ref{find:lotkin}).  We thus showed that a problem admitting no
  solution at all at the fixed precisions of DD/TD/QD can be handled by
  arbitrary-precision \OzII{} in both accuracy and speed.
\item We introduced slice splitting on the CPU as well, choosing the slice
  count $S$ from $S^{*}=\sqrt{A/B}$ of \eqref{eq:TL}.  Without splitting
  ($S=1$) the GEMM grows as $O(p)$ and the conversion as $O(p^2)$, a marked
  imbalance; at the optimum the two grow together towards $p^{1.5}$, making
  the whole $O(p^{1.5})$ (Finding~\ref{find:phase}).  At $n=512$, $p=8192$
  falls from 13.32~s to 3.89~s ($3.4\times$) and $p=16384$ from 31.28~s for
  unsplit binary64 to 10.19~s ($3.1\times$); the $t\lesssim11607$-bit limit
  arising from the INT8 capacity is also removed.  Because slice splitting
  truncates the low-order digit groups, however, problems with severe
  cancellation require a guard sized to the factor $\Gamma$ of
  Proposition~\ref{prop:trunc}.

\item As a consistent trend, higher target precision and larger problem size
  both favour \OzII{}.  This follows directly from the fact that the
  pre-processing and CRT reconstruction costs of \OzII{} do not depend on the
  inner dimension of the Schur complement update, whereas the operation count
  of the low-precision GEMM is proportional to it.
\end{enumerate}

\subsection{Future work: choosing the \OzII{} implementation to match the
number of digits}
\label{sec:future}

The most important lesson of this work is that whether \OzII{} is advantageous
depends on the combination of target precision, problem size and data
representation, and that no single implementation is optimal across the whole
range.  The measured decision points are summarised as follows.

\begin{table}[htbp]
\centering
\caption{Guideline for choosing an implementation, based on measurements
(within the scope of this work)}
\label{tab:selection}
\begin{tabular}{lll}
\toprule
Regime & Recommended & Basis \\
\midrule
DD/TD/QD (CPU, $n\lesssim512$) & machine-dependent (\OzII{} on x86/H100; & Table~\ref{tab:cpu_lu_ddtdqd} \\
 & \quad existing at TD/QD on Arm/GB10) & \\
DD/TD/QD (CPU, $n\gtrsim1024$) & \OzII{} (direct conversion essential) & Finding~\ref{find:cpu_lu_ddtdqd} \\
DD/TD/QD (GPU, weak FP64) & \OzII{}/INT8 & Finding~\ref{find:gpu_reversal} \\
QD (GPU, strong FP64) & \OzII{}/binary64 & Finding~\ref{find:fp64} \\
DD (GPU, strong FP64) & gdtq native & Finding~\ref{find:gpu_reversal} \\
Arbitrary precision (CPU) & \OzII{} (binary64) & Table~\ref{tab:backend_hp} \\
Arbitrary precision (GPU) & \OzII{} (INT8) & Finding~\ref{find:lotkin} \\
FP4 back-end & not recommended (use INT8) & Table~\ref{tab:gpu_gemm} \\
\bottomrule
\end{tabular}
\end{table}

The central task for future work is therefore to establish a mechanism that
automatically selects the \OzII{} implementation appropriate to the number of
digits.  Specifically:

\begin{itemize}
\item A precision- and size-aware selection mechanism: from the target
  precision $p$, the matrix size $n$, the panel width $b$ and the data
  representation (multi-component or arbitrary precision), automatically
  determine the low-precision back-end (INT8/FP16/\dots), the number of slices
  $S$ and the number of moduli $N$, and even whether to use \OzII{} at all,
  based on a cost model and empirical calibration.
  Table~\ref{tab:selection} is a starting point for that calibration data.
\item Accelerating the panel factorisation (highest priority): as
  Finding~\ref{find:profile} shows, after applying \OzII{} the panel
  factorisation accounts for 74--98\% of the total.  Recursive panel
  factorisation, or applying \OzII{} again inside the panel, are possible
  approaches.
\item GPU evaluation at larger $n$: on x86/H100 the whole $n=1024$ LU finishes
  in 0.05~s, making the Schur complement update a small problem of about
  $512^3$.  It is necessary to check how the ranking moves at $n\ge4096$
  (Finding~\ref{find:gpu_reversal}).
\item Implementing the (S0) balancing step, which would make the number of
  moduli $N$ independent of the panel width $b$ and allow $b$ to be chosen
  purely for performance.
\item Choosing the number of slices $S$ for performance: our current CPU
  implementation performs no slice splitting ($S=1$) and obtains the required
  capacity solely by extending the modulus table.  As
  Finding~\ref{find:gemmcount} shows, however, increasing $S$ trades more
  GEMMs for a lower per-element conversion cost, so $S$ is properly a
  performance degree of freedom.  The width of one pass becomes $w=Q/S$, so
  the word counts of both the modular reduction and the CRT fall as $O(1/S)$,
  while the number of GEMMs grows as $O(S)$ by \eqref{eq:gemmcount}.  Writing
  the total time as $T(S)\approx A/S+BL$, the minimum is at
  $S^\ast=\sqrt{A/B}$ with $T=2\sqrt{AB}$.  Substituting the measured values of
  Table~\ref{tab:phase} as $A$ (conversion) and $B$ (GEMM) at $S=1$ gives
  \begin{center}\small
  \begin{tabular}{@{}rrrr@{}}
  \toprule
  $p$ [bit] & $S^\ast$ & $T(S^\ast)$ [s] & vs.\ current \\
  \midrule
  1024 & 3.4 & 0.21 & 1.9 \\
  2048 & 4.3 & 0.47 & 2.3 \\
  4096 & 5.9 & 1.22 & 3.0 \\
  8192 & 7.8 & 3.31 & 4.0 \\
  \bottomrule
  \end{tabular}
  \end{center}
  This is only an estimate from the simple model $A/S+BL$, and it does not
  account for the increase in $Q$ needed to offset the truncation of the lower
  digit groups that accompanies slice splitting
  (Proposition~\ref{prop:trunc}); but given that conversion currently accounts
  for 92--98\% of the total (Finding~\ref{find:phase}), it is well worth
  trying.  Since the GPU implementation already performs slice splitting, this
  amounts not to adding a new mechanism but to porting the GPU structure to
  the CPU and choosing $S$ by performance rather than by the capacity
  condition.

\item Obtaining the \OzII{} breakdown (Table~\ref{tab:phase}) on x86/H100,
  which is missing here because stderr was discarded in that logging run.
\end{itemize}

\section*{Acknowledgements}

This work was supported by JSPS KAKENHI Grant Number JP26K14846.
Anthropic Claude Code and OpenAI ChatGPT were used as assistive tools in
carrying out the implementation and the benchmarks and in preparing this
manuscript.  The author takes full responsibility for the verification of
the experimental results and for the content of this paper.


\end{document}